\documentclass[11pt, a4paper]{amsart}
\usepackage{amsmath, amssymb, amsfonts, amsthm}
\usepackage{color}
\usepackage{bm}
\usepackage{subfigure}
\usepackage{graphicx}
\usepackage{multirow}
\usepackage{setspace}
\usepackage{algorithmic, epsf}
\usepackage{graphicx, color}
\usepackage{paralist}
\usepackage{float}
\usepackage{algorithm}
\usepackage{verbatim}
\usepackage{fullpage}
\usepackage{tikz-cd}
\usepackage{yfonts}
\usepackage{url}
\usepackage[shortlabels]{enumitem}
\usepackage{xpatch}
\title{Instability of the Smith Index Under Joins and Applications to Embeddability}

\author{Salman Parsa}

\address{{University of Utah, Salt Lake City, UT, USA}}
\email{sparsa@sci.utah.edu}
\thanks{{This work was funded by a grant from IPM. This work was funded in part by the National Science Foundation through grant CCF-1614562 as well as funding from the SLU Research Institute.}}

\date{}
\begin{document}

\newcommand{\kfield}        {\mathbf{k}}
\newcommand{\Fspace}        {\mathbb {F}}
\newcommand{\Rspace}        {\mathbb {R}}
\newcommand{\Zring}        {\mathbb {Z}}
\newcommand{\ol}           {\overline}         
\newcommand{\inter}        {\mathcal{X}}
\newcommand{\link}         {\mathcal{L}}
\newcommand{\rf}           {\mathcal{R}^4}
\newcommand{\delembed}           {g^2_\Delta}


\newtheorem{observation}{Observation}
\newtheorem{lemma}{Lemma}
\newtheorem{claim}{Claim}
\newtheorem{Theorem}{Theorem}
\newtheorem{theorem}{Theorem}
\newtheorem*{theorem*}{Theorem}
\newtheorem{proposition}{Proposition}
\newtheorem{corollary}{Corollary}
\newtheorem{Conjecture}{Conjecture}
\newtheorem{question}{Question}

\theoremstyle{definition}
\newtheorem{definition}{Definition}[section]

\theoremstyle{remark}
\newtheorem*{remark}{Remark}

\maketitle

\xpatchcmd{\paragraph}{\normalfont}{{\normalfont\bfseries}}{}{}

\begin{abstract}
    We say a $d$-dimensional simplicial complex embeds into double dimension if it embeds into the Euclidean space of dimension $2d$. For instance, a graph is planar iff it embeds into double dimension. We study the conditions under which the join of two simplicial complexes embeds into double dimension. Quite unexpectedly, we show that there exist complexes which do not embed into double dimension, however their join embeds into the respective double dimension. We further derive conditions, in terms of the van Kampen obstructions of the two complexes, under which the join will not be embeddable into the double dimension. 
    
    Our main tool in this study is the definition of the van Kampen obstruction as a Smith class. We determine the Smith classes of the join of two $\Zring_p$-complexes in terms of the Smith classes of the factors. We show that in general the Smith index is not stable under joins. This allows us to prove our embeddability results. 
\end{abstract}

\section{Introduction}
We say that an $n$-dimensional simplicial complex $M$ \textit{embeds} into (the Euclidean space of) \textit{double dimension} if
there exists a continuous injective map $f:|M| \rightarrow \Rspace^{2n}.$ Conditions for the existence of an embedding into double dimension, for various $n$, is a classical study pioneered by van Kampen and Flores in the 1930's, in attempts to generalize Kuratowski's graph planarity criterion. Our aim in the work presented here is to answer the very basic question:

\begin{question}\label{q:q}
Given two simplicial complexes $M$ and $N$ under which conditions, in terms of $M$ and $N$, does the join simplicial complex $M*N$ embed into double dimension?
\end{question}

An answer to a very special case of this question is given by the well-known classical result which we call the Flores-van Kampen-Gr\"unbaum Theorem \cite{Gru69}, as follows. For $i\in\{1,\ldots,m\}$, let $d_i \geq 0$ be integers and set $d = \sum_i {d_i}+m-1$. Define $K_i$ to be the $d_i$-skeleton of the $(2d_i+2)$-simplex. Then, $$K = K_1 * K_2 * \cdots * K_m,$$ is a $d$-dimensional simplicial complex that does not embed into $\Rspace^{2d}$. The Flores-van Kampen-Gr\"unbaum Theorem is proved by a geometric method and using the Borsuk-Ulam Theorem. We prove a much stronger result by proposing a novel complementary algebraic method based on embedding classes and van Kampen obstructions. 

Our question is more interesting, and has been open, in the case where both $M$ and $N$ do not embed into the Euclidean space of respective double dimensions. We make this assumption for the rest of this introduction.

The answer to Question~\ref{q:q} is known in a few limited cases. The simplest non-trivial case, in a sense, is when one of $M$ or $N$ is 0-dimensional and consists of at least three points. It is known that if the complex $N$ does not embed into the Euclidean space of double dimension, then $[3]*N$ does not embed into double dimension, where by $[3]$ we mean the set of three discrete points. This case can be proved by geometric arguments and appears in the unpublished paper of Melikhov and Shchepin~\cite{MeSh06}. A geometric proof of this fact can be extended to answer the question for all $M$ and $N$ where there is a $\Zring_2$-sphere  
in the deleted join of one of the complexes, see \cite[Theorem~3]{PaSk20} for this extension. We remark that this geometric proof works in the meta-stable dimensions. Unaware of this geometric method, in \cite{Pa20a} the author proved this result by an algebraic method using the theory of Smith classes. Although the algebraic method of \cite{Pa20a} produces a longer proof than the geometric argument in the case of $[3]*N,$ it can be used to answer our question in more generality as presented in this paper.

Another case where the answer is essentially known\footnote{The author was also unaware of this result when publishing \cite{Pa20a}.} is when the modulo 2 van Kampen obstruction classes of $M$ and $N$ are non-zeros. This case has been addressed by Bestvina, Kapovich and Kleiner~\cite[Lemma~9]{BKK02}. Namely, if $M$ and $N$ have non-vanishing van Kampen obstructions modulo 2, then $M*N$ does not embed into the Euclidean space of double dimension\footnote{This result is not explicit in \cite{BKK02}, since the existence of an obstructor appears to be more restrictive than the van Kampen obstruction being non-zero modulo 2. This is because of condition 3 in the definition of an obstructor.}. Another case where the question has been answered is for ``nice'' complexes, see \cite[Exercise 5.8.4]{Mat08} and \cite{Sch93}.

Because of the above results, it is reasonable to conjecture that if the van Kampen obstructions of the factor complexes $M$ and $N$ are non-zeros then the join complex has non-zero obstruction and hence is non-embeddable in double dimension. However, we show in this paper that this impression is not correct. More precisely, we have the following unexpected result.
\begin{theorem}\label{t:existence}
For any $d_M\geq 2$, $d_N \geq 2$ there exist complexes $M$ of dimension $d_M$ and $N$ of dimension $d_N$ such that $M$ does not embed into $2d_M$-dimensional Euclidean space and $N$ does not embed into $2d_N$-dimensional Euclidean space, but $M*N$ embeds into $2(d_M+d_N+1)$-dimensional Euclidean space. $M$ and $N$ can be taken to be any two complexes with non-zero van Kampen obstruction, but zero van Kampen obstructions modulo 2.
\end{theorem}

We also prove that if both $M$ and $N$ have non-zero van~Kampen obstructions modulo 2, then their join does not embed into double dimension\footnote{As we noted, this is not new.}. However, we do not fully answer the question in the case where the obstruction class of one of the complexes vanishes modulo 2 and the other is non-zero modulo 2. In Theorem~\ref{t:mainembed} we give a condition for non-embeddability in this case.

We prove our embeddability results using the concepts of the Smith classes and the Smith index of a $\Zring_2$-complex. These classes give an alternate definition of the van Kampen obstruction, following Wu~\cite{Wu74}. For the purposes of this introduction, it is enough to think of the Smith index, $I(K),$ of a $\Zring_2$-complex $K$, as an integer that is assigned to $K$ and is such that if there is an equivariant map $K \rightarrow S^{m},$ where $S^m$ is the $m$-sphere equipped with the antipodal action, then $I(K)\leq m+1.$ Therefore, the Smith index can serve as an obstruction to embeddability.

In order to prove our results on embeddability of joins, we determine the Smith index\footnote{The Smith index is traditionally defined as the smallest dimension where the Smith class vanishes. For $p=2$ say, this vanishing shows, under some extra conditions, the existence of an equivariant map into an equivariant sphere of dimension one lower than the Smith index. It follows that the Smith index is one greater than the $\Zring_2$-index, as defined in \cite{Mat08} (under those conditions).} of the join of two $\Zring_2$-complexes in terms of the indices of the factors. If $K$ and $L$ are $\Zring_2$-complexes we say that the Smith index is \textit{stable} under the join of $K$ and $L$ if $I(K*L)=I(K)+I(L)$ and \textit{unstable} if the equality does not happen. The surprising embeddability results follow from the fact that the Smith index is unstable in general. We prove our results on the Smith classes for the more general $\Zring_p$-complexes, for $p$ a prime number. Let $I_p(K)$ denote the Smith index of the $\Zring_p$-complex $K$ computed with coefficients $\Zring_p$ (or modulo $p$). Recall that for any $\Zring_p$-complex $K,$ $I(K)-1 \leq I_p(K) \leq I(K).$ Our results on the Smith index of the join can be summarized as follows.

\begin{theorem}\label{t:smith}
Let $K$ and $L$ be $\Zring_p$-complexes and let $I(K)$ and $I(L)$ be their Smith indices. 

\begin{enumerate}
    \item If $I_p(K)=I(K)-1$ and $I_p(L)=I(L)-1$ then $I(K*L) \leq I(K)+I(L)-1$.
    \item If $p=2$, $I_2(K)=I(K)$ and $I_2(L)=I(L)$, then $I_2(K*L)=I(K*L)=I(K)+I(L).$
    \item If $p=2$, $I_2(K)=I(K)-1$ and $I_2(L)=I(L)$ then $I(K*L)=I(K)+I(L)$ or $I(K)+I(L)-1$ and both cases are possible.
    \item If $p>2$,  $I_p(K)=I(K)$ and $I_p(L)=I(L)$, then $I_p(K*L)=I(K*L)=I(K)+I(L)$, given some extra conditions on $K$ or $L$.

\end{enumerate}

\end{theorem}

For details refer to Section~\ref{s:joinsmith} Theorem~\ref{t:stabil}. Previously, the stability of the Smith index was known, for $p=2$, in the case where $M$ or $N$ is a 0-sphere, that is, under suspension \cite{CoFl60}. We have used the term `stability' because of this special case. We expect that analogous theorems to the above theorem can be proved for various homological co-indices defined in \cite{CoFl60}, using our approach.\footnote{After publishing this article, we were notified about the paper~\cite{GoGr2021} which also provides an example were the drop in the cohomological index of the join of two $\Zring_2$-complexes happens.} 

We have also checked our main results by a program that computes the Smith classes for a $\Zring_2$-complex, see the Discussion section for details. 

This paper should be accessible to a large audience. We hope that the more general theorems on the Smith index of $\Zring_p$-complexes, and our method, find application in other problems of combinatorics and discrete geometry whose solutions involve topological methods. For this reason, we have chosen to expose our method with all the details in this paper, rather than to present a short proof of its more interesting outcome, i.e., Theorem~\ref{t:existence}, which has been a by-product of our systematic approach.\\

\paragraph{Method of Proof}
We have developed our method ``from the ground up" building upon simple observations made first for the case of $[3]*K$, as explained in \cite{Pa20a}. This approach has three main components. In the first place, the systematic
definition of the Smith classes in terms of resolutions of the cocycle 1 allows us to calculate formulas for the representatives of Smith classes of the join of two $\Zring_p$-complexes. Second, we are able to show that the cocycles which we have computed belong to non-zero cohomology classes using what we call a ``certificate''. A certificate for a cohomology class is a chain having the property that by evaluating a cocycle of the class on that chain we can deduce that the class is non-zero. We find certificates for Smith classes of the join in terms of the certificates for the Smith classes of the factors. Finally, using deleted joins instead of deleted products simplifies our arguments considerably, and allows us to apply the results on Smith classes to the embeddability question more efficiently. 

It turns out that the concept of a certificate is similar in nature to what is called an obstructor in \cite{BKK02}. However, instead of being combinatorial objects our certificates are algebraic chains, that is, a collection of oriented cells and an integer ``attached" to each one. In a sense, these objects generalize the concept of a minor in a graph. As an example, consider a 3-dimensional complex $M$ and an arbitrary mapping $f:|M| \rightarrow \Rspace^6$ in general position. A certificate for the van Kampen obstruction being non-trivial is a 6-dimensional chain $c$ in the quotient of the deleted product, $c \in C_{6}(\widetilde{M^{\times 2}_\Delta}),$ such 
that the boundary of $c$ is zero modulo $2^m$ for some $m,$ but the intersection homomorphism assigns a number that is not zero modulo $2^m$ to $c$. In any mapping of $M$ into $\Rspace^{6}$ such a chain must exist if the obstruction of $M$ is non-trivial (see Lemma~\ref{l:cert} below). Therefore, in dimensions greater than two, existence of such a chain is equivalent to non-embeddability into double dimension. The certificate chains depend on $f$. However, there are definitions of a cochain in the same cohomology class of the intersection homomorphism that are independent of any map, for instance the Smith class. It follows that 
there are certificates independent of any map. Moreover, if $M$ has a certificate $c$, and $M \subset N$, then $i_\sharp(c)$ also works as a certificate for $N$, where $i$ is the inclusion of deleted products.
There is of course two main differences between minors (or sub-graphs homeomorphic to a minor) and certificates. Namely, the totality of ``minimal" certificates is necessarily infinite \cite{Umm73} (even ignoring the coefficients), and they are chains of the deleted product rather than the complex itself.  

The certificate approach can also be applied in similar situations for demonstrating that a cohomology operation performed on non-zero classes produces a non-zero cohomology class.\\

\paragraph{Sketch of the Proof of Theorem~\ref{t:smith}, Part 1}
Since we deduce Theorem~\ref{t:existence} from this statement, here we present an intuitive justification for this part of Theorem~\ref{t:smith}. The following can serve also as a guide to the actual proof in the body of the paper. Let $K$ and $L$ be $\Zring_2$-complexes, $\Zring_2= \{ 1, t\}$.\footnote{In this paper $\Zring_n = \Zring/n\Zring$.} Denote the homomorphisms induced on chain and cochain groups also by $t$. Let $s=1+t$ and $d=1-t$. The $d$-cochains of a $\Zring_2$-complex are cochains $\phi$ that can be written as $\phi=s\psi$ for some cochain $\psi$. Similarly, the $s$-cochains are those cochains $\phi$ that can be written as $\phi = d \psi$. In both of these cases, let us call $\psi$ a sub-cochain of $\phi$. The Smith classes, $A^i(K)$, of a $\Zring_2$-complex are special cohomology classes represented by a $d$-cocycle in an even dimension and a $s$-cocycle in an odd dimension. Now assume $d_K$ and $d_L$ are highest dimensions where the Smith classes $A_K=A^{d_K}(K)$ of $K$ and $A_L=A^{d_L}(L)$ are non-trivial with integer coefficients, and assume that both are trivial with $\Zring_2$-coefficients. Then $A_K=2B_K$ and $A_L=2B_L$ for some special cohomology classes $B_K,B_L$. For simplicity, assume $d_K$ and $d_L$ are even. It follows from $ss=2s$ that in this situation, we can represent $A_K$ with a cocycle $ss\eta_K$ and $A_L$ with a cocycle $ss\eta_L$. Now our computations show that the tensor product of sub-cochains of the representatives of $A_K$ and $A_L$ is a sub-cochain of a representative of $A(K*L)^{n_K+n_L+1}$. This latter sub-cochain is $(s\eta_K) \otimes (s\eta_L)=s(\eta_K \otimes s\eta_L)$, where we use $s$ also to show the corresponding homomorphism in the cochain complex of join. But since $d_K+d_L+1$ is odd, $A^{d_K+d_L+1}(K*L)$ is represented by $ds(\eta_K \otimes s\eta_L)=0$.\\

\paragraph{Organization of the Paper}
We start in Section~\ref{s:prelim} by recalling some basic notions from the Smith theory which we will need later on in our arguments. In Section~\ref{s:cert}, we prove conditions on chains which can serve as a certificate for Smith classes. In Section~\ref{s:joinresolution}, we relate the Smith classes of the join complex to that of the factors. Section~\ref{s:joinsmith} contains our main theorems on the conditions for the stability and instability of Smith classes. In Section~\ref{s:embed}, we apply our results on Smith classes to the embedding problems and prove our main embedding theorems. Section~\ref{s:rescalc} contains some delayed calculations. At the end we discuss our computer implementation and some open problems.

\section{The Smith Index of a $\Zring_p$-Complex}\label{s:prelim}
In this section, we present the basic background on the concepts we use from the Smith theory of periodic actions of prime period. We restrict ourselves to a brief description of what we shall use in the following sections to prove our results. For a thorough treatment refer to \cite{Wu74,Nak56}. Wu \cite{Wu74} attributes the recasting of the classical results of Smith and Richardson into the cohomology language, which we present in this section, to Thom and Bott.\\

\paragraph{Convention} In this paper, we call what is called a ``special Smith class" in \cite{Wu74} a ``Smith class" and what is called a ``special Smith index" simply ``Smith index". This does not cause any confusion since we do not define nor use the classes that are called ``Smith class" in \cite{Wu74}. In other words, we remove the ``special" adjective throughout, when applied to classes and index.

\subsection{Simplicial $\Zring_p$-Complexes}
Let $K$ be an abstract simplicial complex\footnote{Our complexes are always finite.}  with vertex-set $V(K)=\{v_1, \ldots,v_m \}$. We assume that $K$ is \textit{augmented}, that is, $K$ contains a (-1)-dimensional simplex which we denote by $\emptyset$. As usual $C(K):= \bigoplus_{i=-1}^{\infty} C_i(K)$ is the simplicial chain complex of $K$ with respect to the standard orientation of simplicies of $K$ determined via the ordering of vertices: $v_i<v_j$ iff $i<j$. By definition, for all $i$, $\partial v_i= \emptyset$ and, the homology of $C(K)$ is the reduced homology of $K$.

\begin{definition}
A \textit{simplicial $\Zring_p$-complex} $(K,t)$ is a simplicial complex $K$ together with a bijection $t: V(K) \rightarrow V(K)$ such that i) The simplicial map defined by $t$ permutes the cells of each non-negative dimension, ii) $t^p=id$, where $id$ denotes the identity function. For such $t$ we define $t(\emptyset)=\emptyset$. We denote the simplicial extension of $t$ also by $t:K \rightarrow K$. The action of $t$ is called \textit{simple} if for all non-empty simplices $\sigma$, $t^k\sigma \neq \sigma$ if $k \neq 0 \mod p$.
\end{definition}

\paragraph{Convention} In this paper we work exclusively with simple actions without further notice. Also, we write $K$ instead of $(K,t)$ if no danger of confusion and we abbreviate a simplicial $\Zring_p$-complex to a \textit{$\Zring_p$-complex}. \\


\paragraph{} The prime example of a $\Zring_2$-complex is a simplicial sphere with the antipodal action. The 0-dimensional sphere $S^0$ is just a 2-point set $\{v_1,v_2 \}$ and the empty simplex where $t$ exchanges the two vertices. The $k$-dimensional simplicial sphere $S^k$ with the antipodal action is obtained as the abstract complex resulting from the join $S^0*S^0*\cdots *S^0$ with $k+1$ factors. The permutation $t$ on $V(S^k)$ is simply the product of the permutations of $S^0$ factors. That is, if we denote the vertices of the $j$-th 0-sphere by $\{2j,2j+1\}$, $j \in \{0, \ldots, k \}$, then $t: S^k \rightarrow S^k$ is given by $\Pi_{j=0}^{k}(2j,2j+1).$


\subsection{Special Homologies of a $\Zring_p$-Complex}
The simplicial map $t$ defines a chain homomorphism $t_{\sharp} : C(K) \rightarrow C(K)$ such that $t_\sharp^p=id$. Therefore $C(K)$ is a $\Zring [\Zring_p]$-module, where $\Zring_p= \{1,t_\sharp,\ldots,t_\sharp^{p-1}\}.$ The induced map on the cochain complex is denoted by $t^\sharp: C^*(K) \rightarrow C^*(K)$. $C^*(K)$ is again a $\Zring [\Zring_p]$-module. 

We take two special elements $d,s \in \Zring [\Zring_p]$,

\begin{align*}
     d &= 1-t_\sharp,  \\
     s &= 1+t_\sharp+\cdots + t_\sharp^{p-1}. 
\end{align*}

These two chain homomorphisms have the important property that $d s = s d =0$. For simplicity, we also denote the duals of these homomorphisms which act on the cochain complex $C^*(K)$ by $d$ and $s$. Notice that a cochain in the kernel of $d$ assigns the same value to a chain and its image under $t_\sharp$. More generally, we have the following relations in the chain complex $C(K)$ and also in $C^*(K)$:

\begin{align*}
\ker(d)&= \text{im}(s)\\
\ker(s)&=\text{im}(d).
\end{align*}

See \cite[II.2 Proposition 2]{Wu74} for the proof. If $c \in \ker(d)$ then $d(\partial c) =  \partial (d c)=0$. Hence $\partial c \in \ker(d)$. Analogously, if $c \in \ker(s)$, then $\partial c \in \ker(s)$.

\begin{definition}
Define the chain complex $C^d(K)$ such that $C^d_i(K)=\ker(d_i)$, where $d_i : C_i(K) \rightarrow C_i(K)$ is the restriction of $d$ to $C_i(K)$, and, with the boundary homomorphism which is the restriction of the boundary of $C(K).$  Analogously, $C^s(K)$ is the chain complex such that $C^s_i(K)=\ker(s_i)$ where $s_i: C_i(K) \rightarrow C_i(K)$ is the restriction of $s$ to $C_i(K)$ with the restriction of the boundary homomorphism of $C(K).$

The cochain complex $C^*_d(K)$ is the cochain complex such that $C^i_d(K)=\ker(d_i)$, where $d_i : C^i(K) \rightarrow C^i(K)$ is the restriction of $d$ to $C^i(K)$, and with the same coboundary as $C^*(K)$ but restricted. Analogously, $C^*_s(K)$ is the cochain complex such that $C^i_s(K)=\ker(s_i)$ where $s_i: C^i(K) \rightarrow C^i(K)$ is the restriction of $s$ to $C^i(K)$, with the same coboundary but restricted.

The \textit{$d$-homology}, $H^d(K)$, and the \textit{$s$-homology}, $H^s(K)$, of the $\Zring_p$-complex $(K,t)$ are the homologies of the chain complexes $C^d(K)$ and $C^s(K)$ respectively. The \textit{$d$-cohomology}, $H_d(K)$, and the \textit{$s$-cohomology}, $H_s(K)$, of the $\Zring_p$-complex $(K,t)$ are the (co)homologies of the cochain complexes $C^*_d(K)$ and $C^*_s(K)$ respectively.

\end{definition}

\subsection{The Quotient Complex}

In order to define our main notions, namely the Smith index and Smith classes, it is not necessary to consider the quotient complex. However, these complexes help in relating the $d$-homology and cohomology of $(K,t)$ with an ordinary homology, namely that of the quotient space. This in turn allows to relate the Smith classes to the embedding classes which are traditionally defined in the quotient complex. The quotient complexes also guide our intuition.

\begin{definition}
The \textit{quotient space} of the $\Zring_p$-complex $(K,t)$ is the space $K/t$ which is obtained by identifying each orbit of the action of $t$ on $|K|$ to a single point. If the simplicial structure of $K$ defines a simplicial structure in the quotient space, then the quotient is a simplicial complex and we call it the \textit{quotient complex} of $(K,t)$. 
\end{definition}

If we subdivide the original $\Zring_p$-complex twice then it is guaranteed that the quotient is a simplicial complex (\cite[II.1 Propositions 7 and 8]{Wu74}). We assume this is always the case for our complexes when we consider a quotient space.

The quotient complex $K/t$ has a vertex for every orbit of the $\Zring_p$-action on $V(K)$. We denote this vertex by $[v]\in V(K/t)$, where $v\in V(K)$ is any vertex in the orbit. For simplicity, we assume that the indexing of the $v_i$ is such that the vertices in an orbit are consecutive. Consider a simplex $\sigma = v_0 v_1 \ldots v_m$. The image of this simplex in $K/t$ is the simplex $[\sigma]=[v_0][v_1] \cdots [v_m]$. The simplicial structure of $K/t$ is locally combinatorially equivalent to the structure of $K$. Therefore, 
$$ \partial [\sigma] = \sum_{i=0}^{m} (-1)^i [v_0]\cdots \widehat{[v_i]} \cdots [v_m] = \sum_{i=0}^{m} (-1)^i [v_0\cdots \hat{v}_i \cdots v_m].$$

The natural projection $\pi(v)=[v]$ induces a simplicial projection $\pi: K \rightarrow K/t$. The latter induces in turn the corresponding chain and cochain homomorphisms: $\pi_\sharp : C(K) \rightarrow C(K/t)$ and $\pi^\sharp: C^*(K/t) \rightarrow C^*(K)$. Indeed, one of the goals of studying the $d$- and $s$-homology and cohomology of $(K,t)$ has been to compute the ordinary homology of the quotient space $K/t$. We observe that $\pi_\sharp$ maps any $s$-chain to 0, since a $s$-chain $c$ can be written as $c=(1-t_\sharp)c'$ for some chain $c'$ and $\pi_\sharp((1-t_\sharp)c')=\pi_\sharp(c') - \pi_\sharp t_\sharp(c') = \pi_\sharp(c')- \pi_\sharp(c')=0.$ Moreover, $\pi_\sharp$ cannot define an isomorphism between $d$-homology of $(K,t)$ and homology of $K/t$ since it maps a $d$-chain $sc$ to $p\pi_\sharp (c)$. A more useful map for us is the inverse projection or the transfer homomorphism.

\begin{definition}
Let $\sigma$ be an oriented simplex of $K$, the \textit{inverse projection} $\tau_\sharp: C(K/t) \rightarrow C(K)$ is defined by 

\begin{align*}
\tau_\sharp([\sigma])=  s\sigma = (1+t_\sharp+\cdots+t_\sharp^{p-1})\sigma.
\end{align*}
\end{definition}

Note that $\tau_\sharp$ is well-defined. It assigns to an oriented simplex of $K/t$ the sum of all of its lifts in $K$. The homomorphism $\tau_\sharp$ is a chain homomorphism, \cite[Section 3.G]{Hat02} and \cite[II.1 Proposition~2]{Wu74}.

The projection and the inverse projection and the dual homomorphisms on the cochain complex satisfy the following relations.

\begin{align*}
\tau_\sharp \pi_\sharp &= s,\;\; \pi_\sharp \tau_\sharp = p\cdot id,\\
\pi^\sharp \tau^\sharp &= s,\;\;  \tau^\sharp \pi^\sharp = p \cdot id.
\end{align*}

Observe that $\tau_\sharp(C(K/t)) \subset C^d(K)$. In fact, $\tau_\sharp(C(K/t))=C^d(K).$ Let $\bar{\tau}_\sharp$ be equal to $\tau_\sharp$ but with the co-domain restricted to $C^d(K).$ $\bar{\tau}_*: H(K/t) \rightarrow H^d(K)$ is an isomorphism (\cite[II.2 Proposition 4]{Wu74}). The dual homomorphism $\bar{\tau}^\sharp: \text{Hom}(C^d(K), \Zring) \rightarrow C^*(K/t)$ is a chain homomorphism and is an isomorphism on each cochain group. Note that $\bar{\tau}^\sharp$ is not an isomorphism of $C^*_d(K)$ and $C^*(K/t)$\footnote{ And therefore in \cite[II.2 proposition 4]{Wu74} $\tilde{\tau}^\sharp$ has to be substituted for $\bar{\pi}^\sharp$.}. Hence we need yet another homomorphism $\tilde{\tau}^\sharp: C^*_d(K) \rightarrow C^*(K/t)$. Let $\phi \in C^*_d(K)$ and $\psi$ be such that $\phi = s \psi$. Then define

$$ \tilde{\tau}^\sharp(\phi) = \tau^\sharp(\psi).$$

Further, let $\tilde{\pi}^\sharp: C^*(K/t) \rightarrow C^*_d(K)$ be the homomorphism induced by $\pi$ and with the co-domain restricted.

\begin{lemma}\label{l:dquotientiso}
The function $\tilde{\tau}^\sharp$ is a chain homomorphism. In addition, $\tilde{\tau}^*$ and $\tilde{\pi}^*$ are inverse isomorphisms of $H^*_d(K)$ and $H^*(K/t).$

\end{lemma}

\begin{proof}
First we show that $\tilde{\tau}^\sharp$ is well-defined. Let $\psi'$ be another cochain such that
$\phi = s \psi'$. Then $s(\psi-\psi')=0$, therefore $\psi-\psi' = d \eta$ for some cochain $\eta$ and $\tau^\sharp(d \eta)=0$, and thus $\tilde{\tau}^\sharp(\psi)=\tilde{\tau}^\sharp(\psi').$ Also one checks that $\tilde{\tau}^\sharp$ is a homomorphism.

We next show that $\tilde{\tau}$ is a chain homomorphism. We have $\delta (\tilde{\tau}^\sharp (s\psi)) = \delta( \tau^\sharp (\psi)) = \tau^\sharp (\delta (\psi)) = \tilde{\tau}^\sharp (s \delta (\psi)) = \tilde{\tau}^\sharp (\delta (s\psi)).$ To see that $\tilde{\tau}^\sharp$ is an isomorphism one checks that $\tilde{\pi}^\sharp \tilde{\tau}^\sharp = id_{C^*_d(K)}$ and $\tilde{\tau}^\sharp \tilde{\pi}^\sharp = id.$
\end{proof}

\begin{figure}
    \centering

    \begin{tikzcd}
    C^d(K) \arrow[r,hook] & C(K) \arrow[d, "\pi_\sharp"] & \text{Hom}(C^d(K),\Zring) \arrow[rd, "\bar{\tau}^\sharp\;\cong"']  & \arrow[l] C^*(K) \arrow[d,shift left=1ex, "\tau^\sharp", end anchor={[xshift=1ex]}, start anchor={[xshift=1ex]}] & \arrow[l,hook'] C^*_d(K) \arrow[ld, bend left, "\tilde{\tau}^\sharp\;\cong", end anchor={[xshift=1ex]}, start anchor={[xshift=1ex]}]\\
    & C(K/t) \arrow[ul,"\bar{\tau}_\sharp\; \cong"] \arrow[u,"\tau_\sharp",shift left =1.5ex] & & C^*(K/t) \arrow[u, "\pi^\sharp"] \arrow[ur, bend right= 20, "\tilde{\pi}^\sharp\; \cong"] 
    \end{tikzcd}
    \caption{Homomorphisms of this section}
\end{figure}

\paragraph{}
The $s$-homology and $s$-cohomology of a $\Zring_p$-complex can also be described as a homology of the quotient complex $K/t$. Let $\Zring^p=\bigoplus_{i=1}^{p} \Zring.$ We denote by $\Zring^p_0 \subset \Zring^p$ the subgroup of all elements $(n_1, \ldots, n_p) \in \Zring^p$ such that $\sum_{i=1}^{p} n_i = 0.$ Then $H^s(K)$ is isomorphic with \textit{a} homology of $C(K/t)$ with coefficients in $\Zring^p_0$ and analogously $H^*_s(K)$ is isomorphic with \textit{a} cohomology of $C(K/t)$ with coefficients in $\Zring^p_0$, for more details refer to \cite[II.2 Proposition 5]{Wu74}.

\subsection{The Smith Classes of a $\Zring_p$-Complex}\label{ss:smithclass}
Let $\phi_0 \in C_d^i(K)$ be a $d$-cocycle. Then there is $\psi_0$ such that $\phi_0 = s \psi_0$. We have
$0=\delta \phi_0 = \delta s \psi_0 = s \delta \psi_0 = 0.$ Therefore, $\delta \psi_0$ is an $s$-cocycle.
We can write $\delta \psi_0= d \psi_1$ for some cochain $\psi_1$. Then again $0=\delta\delta \psi_0= \delta d \psi_1 = d \delta \psi_1.$ Thus $\delta \psi_1$ is a $d$-cocycle and we can write $\delta \psi_1 = s \psi_2$. We can continue this process indefinitely. We could have started with an $s$-cocycle. The sequence of the $\psi_j$ obtained in this way is called a resolution for the cocycle $\phi_0$.

\begin{definition}
Let $\phi \in C_d^i(K)$ be a $d$-cocycle. A \textit{resolution} for $\phi$ is a sequence of cochains $\psi_j \in C^{i+j}(K), j=0,1,2,\ldots$ such that i) $\phi = s \psi_0$, ii) if $j$ is even,  $\delta \psi_j = d \psi_{j+1}$ and iii) if $j$ is odd, $\delta \psi_j = s \psi_{j+1}.$ A resolution for an $s$-cocycle is defined analogously.
\end{definition}

Since $\delta \psi_j$ is an $s$- or $d$-cocycle we can take its cohomology class in the corresponding $s$- or $d$-cohomology. The resulting cohomology classes are independent of the chosen resolution and the starting cocycle. Thus a resolution can be used to define homomorphisms that push a $d$- or an $s$-cohomology class up to higher dimensions and alternatively to $d$- and $s$- cohomologies. We do not make use of these \textit{special Smith homomorphisms} in this paper and refer the interested reader to \cite[Chapter~II]{Wu74}.

Let $1 \in C^0(K)$ be the cochain that assigns the value $1 \in \Zring$ to each vertex. Then $1$ is a $d$-cocycle.

\begin{definition}
Let $\psi_j \in C^j(K), j=0,1,\ldots$ be a resolution for $1 \in C^0(K)$. \textit{The $j$-th
Smith class} of $(K,t)$, denoted $A^j(K,t)$, is the $d$-cohomology class of $s \psi_{j}$ if $j$ is
even and is the $s$-cohomology class of $d \psi_{j}$ if $j$ is odd. Therefore, for $j$ even $A^j(K,t) \in H^*_d(K)$ and for $j$ odd $A^j(K,t) \in H^*_s(K).$

The \textit{Smith index} of a $\Zring_p$-complex $(K,t)$ is the smallest $n$ such that
$A^n(K,t)=~0.$ We denote the Smith index by $I(K,t)$ and by $I(K)$ when the action is not important.
We call the non-trivial Smith class $A^{I(K,t)-1}(K,t)$ the \textit{top Smith class}. We denote by $A_p^{j}(K,t)$ the Smith classes computed modulo a number $p$,\footnote{That is, $A^j_p(K)$ is the image of $A^j(K)$ under reduction of the coefficient group modulo $p$.} and by $I_p(K,t)$
the smallest integer such that $A_p^{j}(K,t)=0.$
\end{definition}

\paragraph{Remark}
Since our complexes are augmented, we could define the resolution of 1 starting from dimension -1, by setting $\psi_{-1}= \emptyset.$ Then $\delta(\emptyset)=1=s\psi_0.$

\begin{proposition}
Let $(K,t), (K',t')$ be two $\Zring_p$-complexes and let $f:K \rightarrow K'$ be an equivariant simplicial map. Then $f$ induces chain and cochain maps on $s$- and $d$-chain and cochain complexes and we have $$f_j^*(A^j(K',t'))=A^j(K,t)$$
where $f_j^*:H^*_d(K') \rightarrow H^*_d(K)$ if $j$ is even and $f_j^*:H^*_s(K') \rightarrow H^*_s(K)$ if $j$ is odd. It follows that if such an $f$ exists $I(K,t) \leq I(K',t').$ 
\end{proposition}

As an example for $p=2$ we can take any $\Zring_2$-complex $(K,t)$ and map it via a simplicial map to a sufficiently high dimensional sphere with the antipodal action. The resulting map is unique up to equivariant homotopy. Therefore the Smith classes are the images of the Smith classes of the infinite-dimensional sphere via a unique map on $d$- and $s$-cohomologies. Therefore, these classes submit to the same additional natural structure as the classes of the infinite-dimensional sphere and the corresponding classes of the quotient complex follow the structure of the classes of the infinite-dimensional projective space.
For instance, $\tilde{\tau}^\sharp A^{2k} \smallsmile \tilde{\tau}^\sharp A^{2l} = \tilde{\tau}^\sharp A^{2k+2l}.$ See \cite[Section II.2]{Wu74} and \cite{Sha57}. Analogously for other $p>2$ one can take a classifying space and deduce the structure of the Smith classes. However the $s$-cohomologies introduce many complications. For this reason they are sometimes mapped into $\Zring_p$-homology of the quotient complex, as in \cite{Wu74}.  

The Smith classes of positive dimension are $p$-torsion elements of the corresponding cohomology groups.
\begin{lemma}
    For any $i \geq 1$ the Smith classes satisfy $pA^i(K,t)=0.$
\end{lemma}

\begin{proof}
Let $\psi_i$ be a resolution for $1 \in C^*_d(K).$
Let $i=2k$ be even. Recall that $A^{2k} = [\delta \psi_{2k-1}]_d.$ where $[\cdot]_d$ denotes the $d$-cohomology class. Also $\delta (\psi_{2k-1}) = s \psi_{2k}.$ Since $ss=ps$ we have $\delta (s \psi_{2k-1}) =  s \delta (\psi_{2k-1})= s\cdot s \psi_{2k} = p s \psi_{2k} = p \delta (\psi_{2k-1})$. Therefore, $p \delta (\psi_{2k-1})$ is a $d$-coboundary.

Next assume $i=2k+1$ is odd. Then $A^{2k+1} = [\delta \psi_{2k}]_s.$ Also $\delta \psi_{2k} = d \psi_{2k+1}.$
We have the relation $(p-s)(1-t^\sharp)=p(1-t^\sharp)$ and in our notation $(p-s)d = pd.$ Therefore,
$\delta ((p-s)\psi_{2k}) = (p-s)d \psi_{2k+1} = pd \psi_{2k+1}=p \delta \psi_{2k}.$ Moreover, a calculation shows that $p-s = (1-t^\sharp)(p-1 + (p-2)t^\sharp + \cdots + 2 (t^\sharp)^{p-3} + {t^\sharp}^{p-2}).$ Thus
$(p-s) \psi_{2k}$ is an $s$-cochain and it follows that $p[\delta (\psi_{2k})]_s=0.$

\end{proof}

\section{Certificate for a Non-trivial Cohomology Class }\label{s:cert}
In order to prove our main results we need to show that the cohomology classes of certain cocycles are non-zero. Our approach is to demonstrate a chain such that from the evaluation of the cocycle on that chain we can deduce that the cohomology class of the cocycle is non-trivial. This is easy when the coefficients are from a field (Kronecker pairing) or when the cohomology class we are interested in is not a torsion element. However, our cohomology classes are torsion elements. Such elements assign zero to all the cycles. Nevertheless, it is still possible to certify that the class of a cochain is a non-trivial torsion cohomology class by evaluation of the cocycle on a chain. The following lemma seems to be part of homological algebra. However since we have not encountered it in the literature in this form we provide the details of the argument. 

\begin{lemma}\label{l:cert}
    Let $C$ be a chain complex of free abelian groups with integer coefficients, having a finite set of generators in each dimension. Let $h \in H^i(C)$ be a torsion cohomology class. Then, $h \neq 0$ if and only if there exist an integer $n \geq 2$ and a chain $c \in C_i$ such that $c$ is a cycle modulo $n$ and $\phi(c)\neq 0$ modulo $n$ for some representative $\phi$ of $h$. Moreover, if for some prime $p$ and $m\geq 1$, the reduction modulo $p^m$ of $h$ is non-zero one can take $n=p^m$.
\end{lemma}

\begin{proof}
If such an integer $n$ and chain $c$ exist then $h\neq 0$ in cohomology modulo $n$ hence also as an integer. We thus focus on the other direction.
Let $D:C_{i} \rightarrow C_{i-1}$ denote the $i$-dimensional boundary matrix of $C$. We first write the Smith normal form: $D=VSU$ where $U$ and $V$ are non-singular integer matrices with determinant $\pm 1$. Let $u_i$ denote the columns of $U$. $S$ is the matrix of the boundary operator with the new bases $\{u_i\}$ for $C_i$ and the basis consisted of columns of $V^{-1}$ for $C_{i-1}$. Let $\eta_i \in C^i$ be the cochain that assigns the value 1 to $u_i$ and 0 to other basis elements. Similarly let the $\theta_j \in C^{i-1}$ form the basis dual to the basis defined by $V^{-1}$. Then, the coboundary matrix $\delta: C^{i-1} \rightarrow C^{i}$ equals $S^t$ with respect to the dual bases $\{\eta_i\},\{\theta_j\}$. Let $S^t = \text{diag}(1,\cdots,1,s_a,s_{a+1},s_{a+2},\ldots,s_b,0,0,\ldots)$, where $s_i|s_{i+1}$ and $s_a>1$.  Note that this implies $\delta \theta_j = s_j\eta_j$, where $s_j\neq0$, hence $\delta s_j\eta_j=0$ and thus $\delta\eta_j=0$ for those $j$. It follows that $\{[\eta_{i}], a\leq i \leq b  \}$ generate the torsion subgroup of $H^i(C)$ and order of $[\eta_{i}]$ equals $s_i$. 

We choose the $u_i$ to be a certificate for $k[\eta_i]$, $0<k<s_i, a \leq i \leq b$. Note that $\partial u_i= s_i x=0 \mod s_i$ for some chain $x$. Moreover, $\eta(k u_i)=k\eta(u_i)=k \neq 0 \mod s_i$.

Now let $h$ be an arbitrary torsion cohomology class. Then, $h= \sum_{j=a}^b m_j [\eta_j]$ where $0 \leq m_j < s_j$. It is easily seen that $n=s_j$ and $c=u_j$, for any $j$ such that $m_j \neq 0$, is a certificate for $h$.   

For any given $h$ we have demonstrated $c$ and $n$ for some $\phi$ in the statement of the lemma. Let $\phi'$ be another representative for $h$. Then $\phi' - \phi = \delta \psi$ for some cochain $\psi$ and $\phi'(c)=\phi(c) + \psi( \partial c)$. It follows that $\phi'(c)=\phi(c) \mod n$. 

We now prove the last statement of the lemma. Let $h = \sum_i {n_i [\eta_i]}$. Recall that the number $s_{a+i}, 0 \leq i \leq b-a$, is the product of elementary divisors $p_j^{m_j}$ such that, for each $j$, there are at least $b-i$ other elementary divisors larger than or equal to $p_j^{m_j}$ and that for each prime $p_j$ there is at most a single factor $p_j^{m_j}$ in $s_{a+i}$. There is some $k$ such that $n_k \eta_k$ is non-trivial modulo $p^m$. Note that $\eta_k$ is a $s_k$-torsion and since it is not trivial modulo $p^m$ it must be that $s_k$ has a factor $p^{m'}$ for some $m'$ and $n_i\neq 0 \mod p^m$. We observe that, in addition, $n_k \neq 0 \mod p^{m'}$, since $n_k \eta_k$ is mapped to $n_k(1 \oplus 1 \oplus \cdots \oplus 1)=\bigoplus_j (n_k \mod p_j^{m_j})$, under the isomorphism $$\Zring/s_k\Zring \cong \bigoplus_j \Zring/p_j^{m_j}\Zring,$$ and $p^{m'}=p_{j_{0}}^{m_{j_0}}$ for a single $j_0$.

Assume $m \geq m'$. Consider the certificate $u_k$ of $\eta_k$ and note that $\partial (p^{m-m'}u_k) = p^{m-m'} s_kx$. It follows that $\partial p^{m-m'} u_k = 0 \mod p^{m}$. Moreover $n_k \eta_k (p^{m-m'}u_k) = n_k p^{m-m'} \neq 0$ all modulo $p^m$. The last relation follows since $n_k\neq 0 \mod p^{m'}$. Therefore $p^{m-m'} u_k$ is a certificate for this case with $n=p^m$.

Assume $m < m'$. Then $\partial (u_k) = s_kx=0 \mod p^m$. Moreover, $h(u_k)=n_k \neq 0 \mod p^m$. We have a certificate in this case too.

\end{proof}

\paragraph{Examples} To illustrate the preceding proposition we give some examples. 

1) As a simple example, consider the real projective space $\Rspace P^2$ and the cell complex consisting a single cell $e^i$ in dimension $i$. The co-chain $\phi$ that assigns 1 to $e^2$ generates the cohomology group $H^2(\Rspace P^2)= \Zring_2$. We can certify that $[\phi]\neq 0$ using $e^2$. We have $\partial e_2 = 2 e_1=0 \mod 2$ and $\phi(e_2)\neq 0 \mod 2$. 

2) As a second example, let $X$ be a space similar to $\Rspace P^2$ but the the disk $e^2$ is now attached using a degree-4 map instead of a degree-2 one. Then $H^2(X)=\Zring_4$, generated by the class of the cochain $\phi$ that assign 1 to $e^2$. Let's find a certificate for $[\phi]$ modulo $2^3=8$. Consider the chain $2e$. We have $\partial 2e= 8e_1=0 \mod 8$. On the other hand, $\phi(2e)=2 \neq 0 \mod 8$. Thus $2e$ is a certificate for $[\phi]$. The same chain is also a certificate for the classes $2[\phi]$ and $3[\phi]$ modulo 8.

3) As a third example, take $X$ but we now we want to show that $3[\phi]$ is not trivial modulo 2. We take the chain $e_2$ as a certificate. Then $\partial(e_2)=4e_1=0 \mod 2$ and $3\phi(e_2)=3 \neq 0 \mod 2$. Note that when reduced modulo 2, $3\phi$ and $\phi$ belong to the same class.

\begin{proposition}
For a chain complex $C$ as in Lemma~\ref{l:cert}, If $h\neq 0$ is a $p$-torsion cohomology class then there is an $m$ such that $h \neq 0$ modulo $p^m$.
\end{proposition}
\begin{proof}
From the Bockstein long exact sequence it follows that if $h$ is non-trivial and the reduction modulo $p$ is trivial, then $h$ can be written as $h=ph_1$ for some $h_1 \neq 0.$ Similarly for $h_1$, if it is not zero modulo $p$ it can be written as $h_1=ph_2.$ We can continue this process such that $h=p^kh_k.$ The $h_k$ are torsion elements and there are a finite number of them. Therefore it follows that for some $k_1>k_2$, $h_{k_1}=h_{k_2}$ and $h = p^{k_1}h_{k_1}=p^{k_1}h_{k_2}=p^{k_1-k_2}p^{k_2}h_{k_2}=p^{k_1-k_2}h=0$, which is a contradiction. Thus for some $m'$, $h_{m'} \neq 0 \mod p$. Then $h = p^{m'}h_{m'}\neq 0 \mod p^{m'+1}$ since otherwise $h=p^{m'+1}h'=p^{m'}h_{m'}$ for some nonzero $h'$, and $ph'=h_{m'}$ which is a contradiction. 

\end{proof}

\subsection{Certificate for a Non-trivial $d$-Cohomology Class}
Our next goal is to find certificates for non-trivial $d$- and $s$-cohomology classes of $(K,t)$ in terms of chains of $C(K).$ Our main task here is to express the $d$-cohomology and $s$-cohomology as an ordinary cohomology of a free chain complex and then use Lemma~\ref{l:cert}.

We start by defining an isomorphism $C^*_d(K) \cong \text{Hom}(C^d(K), \Zring).$ Let $\{ \sigma_j, j \in J \}$ be a specific set of oriented simplices that form a fundamental domain for all the simplices of $K$, where $J$ is an index set. 
Then $\{t^k_\sharp {\sigma}_j, j \in J, k=0,\ldots, p-1\}$ is a set of free generators for $C(K)$ (in non-negative dimensions). We briefly denote by $a^*$ the dual cochain to a free generator $a,$ that is, the cochain that assigns 1 to $a$ and 0 to other generators. The set $\{s \sigma^*_j, j \in J\}$ is a set of free generators for $C^*_d(K)$ and the set  $\{s {\sigma}_j, j \in J \}$ is a set of free generators for $C^d(K).$
Define $\gamma=\gamma_d: C^*_d(K) \rightarrow\text{Hom}(C^d(K), \Zring)$ on $s\sigma_j^* \in C^*_d(K)$ as $$\gamma (s\sigma_j^*) = (s\sigma_j)^*.$$ 

\begin{lemma}\label{l:gammahom}
The homomorphism $\gamma: C^*_d(K) \rightarrow\text{Hom}(C^d(K), \Zring)$ is an isomorphism of chain complexes. Moreover, for any $s \phi \in C^*_d(K)$ such that support of $\phi$ is in the fundamental domain, and any chain $sc \in C^d(K),$ $\gamma(s\phi)(sc)=\phi(sc).$
\end{lemma}

\begin{proof}
It is clear that $\gamma$ is an isomorphism on each cochain group. Moreover, it is straightforward to check that $\gamma = (\bar{\tau}^\sharp)^{-1} \tilde{\tau}^\sharp,$ therefore it is a chain homomorphism. To prove the last statement, write $\phi = \sum_j n_j  \sigma_j^*$ and $sc= \sum_i m_i s {\sigma_i}.$ We compute,
\begin{align*}
    \gamma(s\phi)(sc) &= \sum_j n_j \left(s{\sigma}_j\right)^*\left(\sum_{i} m_i s{\sigma_i}\right)\\
    &= \sum_j n_jm_j.
\end{align*}

On the other hand,
\begin{align*}
    \phi(sc) &= \sum n_j \sigma_j^* \left(\sum m_i s {\sigma_i}\right)= \sum_j n_jm_j.
\end{align*}
\end{proof}

The following describes the conditions on a chain in $C(K)$ to serve as a certificate for a non-trivial $d$-cohomology class.

\begin{lemma}\label{l:certd}
    Let $C(K)$ be the chain complex of a $\Zring_p$-complex $K.$ Let $h \in H_d^i(K)$ be a torsion $d$-cohomology class and assume $h=[s \phi]_d$, where the support of $\phi$ is in a fundamental domain. Then, $h \neq 0$ if and only if there exist an integer $n \geq 2$ and a chain $c \in C_i(K)$ such that,
    
    \begin{enumerate}
        \item[i)] $c$ is a linear combination of oriented simplices in the fundamental domain containing the support of $\phi$,
        \item[ii)] $s \partial(c) =0$ modulo $n,$
        \item[iii)] $\phi(sc)=\phi(c)\neq 0$ modulo $n$.
    \end{enumerate}

    Moreover, if for some prime $p$ and $m\geq 1$, the reduction modulo $p^m$ of $h$ is non-zero one can take $n=p^m$.
\end{lemma}

\begin{proof}
Assume such $n$ and $c$ exist. Note that $\partial (sc) = s\partial (c)=0$ modulo $n$ and $sc \in C^d(K).$ Moreover, $\gamma (s\phi) (sc) = \phi(c) \neq 0$ modulo $n$. Therefore $[\gamma (s\phi)]$ is non-trivial by Lemma~\ref{l:cert} and so is the $d$-cohomology class of $s\phi.$

To prove the other direction let $n$ and $c'$ be certificates for $[\gamma(s\phi)]$ given by Lemma~\ref{l:cert}. Then $c'$ is a $d$-chain and can be written as $c' = s c$ for some $c \in C(K)$ where $c$ is a linear combination of oriented simplices in the fundamental domain. Then $ \phi(c) = \gamma (s\phi) (sc) = \gamma (s\phi) (c') \neq 0$ modulo $n$. Moreover, $s \partial (c) = \partial (sc) = \partial c' = 0$ modulo~$n$.
The last statement follows from the last statement of Lemma~\ref{l:cert}.

\end{proof}

\subsection{Certificate for a Non-trivial $s$-Cohomology Class}

We proceed to find certificates for $s$-cohomology classes among chains of $C(K).$ There is an exact sequence of chain complexes,
$$0 \rightarrow C^s(K) \xrightarrow[]{\iota_s} C(K) \xrightarrow[]{j_ds} C^d(K) \rightarrow 0,$$
where $\iota_s$ denotes the inclusion and $j_d$ is a projection of $C(K)$ onto $C^d(K)$ which is the identity on $C^d(K)$. Since $C^d(K)=\ker(d)$ is a free subgroup the above sequence is split (in each dimension) and we obtain a short exact sequence of cochain complexes

$$0 \rightarrow \text{Hom}(C^d(K),\Zring) \xrightarrow[]{sj_d^\sharp} C^*(K) \xrightarrow[]{\iota^\sharp_s} \text{Hom}(C^s(K), \Zring) \rightarrow 0.$$
Note that here we have used the split property of the exact sequence in each dimension to show that
$\text{Hom}(,)$ applied in each dimension produces a split short exact sequence, see, for instance, \cite{Spa89}[Chapter 5, Lemma 4.7].

Define the homomorphism $\gamma_s : C^*_s(K) \rightarrow \text{Hom}(C^s(K), \Zring)$ by the formula 
$$ \gamma_s (d \psi) = \iota^\sharp_s (\psi), \; \; d\psi \in C^*_s(K), \psi \in C^*(K).$$
Also define $\rho_s : \text{Hom}(C^s(K),\Zring) \rightarrow C^*_s(K)$ by the formula
$$ \rho_s (\iota^\sharp_s(\eta)) = d \eta, \; \; \iota^\sharp_s (\eta) \in \text{Hom}(C^s(K),\Zring), \eta \in C^*(K).$$

\begin{lemma}
    The homomorphisms $\gamma_s$ and $\rho_s$ are well-defined chain homomorphisms and they induce
    isomorphisms on homology $H^*(C^s(K)) \cong H^*_s(K)$.
\end{lemma}

\begin{proof}
We first show that $Im(sj_d^\sharp)= C^*_d(K).$ Cleary $Im(sj_d^\sharp) \subset C^*_d(K)=s(C^*(K)).$ On the other hand, if $x \in C^s(K),$ $\iota^\sharp_s(s\phi)(x) = s\phi (\iota_s(x)) = s\phi(dy) =0$. Thus $C^*_d(K) \subset \ker(\iota^\sharp)=Im(sj_d^\sharp).$

Consider $\gamma_s.$ If $d \psi = d \psi'$ then $d (\psi - \psi')=0$ thus $\psi - \psi' = s \phi$ for some $\phi \in C^*(K).$ It follows that $\gamma_s(d \psi)  -  \gamma_s(d \psi' ) = \iota_s^\sharp (\psi-\psi') = \iota^\sharp_s (s \phi ) =0.$ Next consider $\rho_s.$ If $\iota^\sharp_s \eta = \iota^\sharp_s \eta'$ then $\eta - \eta' = sj_d^\sharp (\zeta)$ for some $\zeta \in \text{Hom}(C^d(K),\Zring).$ It follows that
$\rho_s(\iota^\sharp_s \eta) - \rho_s(\iota^\sharp_s \eta') = d sj_d^\sharp (\zeta)=0.$ Therefore these two homomorphisms are well-defined.

We calculate $\delta \gamma_s (d \psi)=  \delta \iota^\sharp_s (\psi) = \iota^\sharp_s \delta (\psi)
= \gamma_s ( d \delta (\psi)) = \gamma_s \delta (d\psi)$. Similarly $\delta \rho_s(\iota^\sharp_s(\eta))=
\delta d \eta = d \delta \eta = \rho_s(\iota^\sharp (\delta \eta)) = \rho_s \delta (\iota^\sharp_s(\eta)).$ 

Using the homomorphisms $\rho_d = \gamma_d^{-1}$ of the previous section we define a homomorphism of exact sequences 
\[
  \begin{tikzcd}
    0 \arrow{r} & \text{Hom}(C^d(K),\Zring) \arrow{d}{\rho_d} \arrow{r}{sj^\sharp_d} & C^*(K) \arrow[d,equal]{} \arrow{r}{\iota^\sharp_s} & \text{Hom}(C^s(K),\Zring) \arrow{d}{\rho_s}  \arrow{r}{} & 0 \\
     0 \arrow{r}{} & {C^*_d(K)} \arrow[r,hook]{}  & C^*(K) \arrow{r}{d} & C^*_s(K) \arrow{r}{} &0.
  \end{tikzcd}
\]

For every $ \eta \in \text{Hom}(C^d(K),\Zring)$ and every oriented simplex $\sigma$ we have $sj^\sharp_d(\eta)(\sigma)= \eta(j_d s \sigma) = \eta(s\sigma) = \rho_d(\eta)(\sigma).$ Therefore the left square in the diagram is commutative and so is the second one. Since $\rho_d$ induces isomorphisms on (co)homology groups, the Five Lemma applied to the induced homomorphisms of the long exact sequences of the above short exact sequence of chain complexes shows that $\rho_s$ also induces isomorphisms on (co)homology groups. Since $\rho_s \gamma_s = id$ and $\gamma_s \rho_s =id$, $\gamma_s$ also induces isomorphisms on (co)homology groups.
\end{proof}

We are now ready for describing the chains that serve as certificates for non-trivial $s$-cohomology classes.

\begin{lemma}\label{l:certs}
    Let $C(K)$ be the chain complex of a $\Zring_p$-complex $K.$ Let $h \in H_s^i(K)$ be a torsion $s$-cohomology class and assume $h=[d \phi]_s$. Then, $h \neq 0$ if and only if there exist an integer $n \geq 2$ and a chain $c \in C_i(K)$ such that $d \partial(c) =0$ modulo $n$ and $\phi(dc)\neq 0$ modulo $n$. Moreover, if for some prime $p$ and $m\geq 1$, the reduction modulo $p^m$ of $h$ is non-zero one can take $n=p^m$.
\end{lemma}

\begin{proof}
Assume such $n$ and $c$ exist. Note that $\partial (dc) = d\partial (c)=0$ modulo $n.$ Since $\gamma_s (d\phi) = \iota^\sharp_s \phi$, it follows that $\gamma_s (d\phi)(dc)= (\iota^\sharp_s \phi)(dc)=\phi (\iota_s dc)=\phi(dc) \neq 0 \mod n$. It follows that $[\gamma_s (d\phi)] \neq 0$ and since $\gamma_s$ induced isomorphism on homology $[d\phi]_s$ is non-trivial as an $s$-cohomology class.

To prove the other direction let
$n$ and $c'$ be given by Lemma~\ref{l:cert} for $[\gamma_s(d\phi)]$. Then $c'$ is an $s$-chain and can be written as $c' = d c$ for some $c \in C(K)$. Then $\phi(dc)= \phi(\iota_s dc)= (\iota^\sharp_s \phi) (dc) = (\gamma_s d \phi) (dc) \neq 0$ modulo $n$. Moreover, $d \partial c = \partial dc = \partial c' = 0$ modulo $n$. 

The last statement follows from the last statement of Lemma~\ref{l:cert}.
\end{proof}

\section{Joins of $\Zring_p$-Complexes and Their Resolutions}\label{s:joinresolution}

Let $(K,t_K)$ and $(L,t_L)$ be $\Zring_p$-complexes. The join $K*L$ is a $\Zring_p$-complex
with the action given by $ t: V(K*L) \rightarrow V(K*L)$,  $t(v)=t_K(v)$ if $v \in V(K)$ and $t(v)=t_L(v)$ if 
$v \in V(L)$ where we assume $V(K*L)=V(K) \sqcup V(L).$ Of course, when considering complexes like $K*K$ we rename the vertices of the second copy of $K.$ When orienting the simplices of $K*L$ we use the orderings of the vertices used to define $C(K)$ and $C(L)$ but place vertices of $L$ after those of $K.$ In this way, the orientations of simplices of $K$ and $L$ define an orientation on simplices of $K*L.$ Our goal in this section is to compute a resolution of $1 \in C^*_d(K*L)$ in terms of resolutions in $K$ and in $L.$ \\

\paragraph{}
Recall that our simplicial complexes are augmented. Define the chain complex $D(K*L)$ by
$$ D_n(K*L) = \left(C(K)\otimes C(L) \right)_{n-1} = \sum_{j} C_{n-1-j}(K) \otimes C_{j}(L),$$
and the boundary operator given by
$$\partial(\sigma \otimes \tau) = \partial \sigma \otimes \tau + (-1)^{d_\sigma+1}\sigma\otimes \partial \tau,$$
where $d_\sigma$ denotes dimension of $\sigma.$ The chain complex $D$ is also used by Milnor~\cite{Mil56}.

\begin{proposition}
Let $C(K)$ and $C(L)$ be chain complexes of augmented simplicial complexes $K$ and $L$. Let $\sigma \in C(K)$ and $\tau \in C(L)$ and $\sigma*\tau \in C(K*L)$ be oriented simplices. Then the homomorphisms induced by  $\sigma*\tau \rightarrow \sigma \otimes \tau$ and $\sigma \otimes \tau \rightarrow \sigma * \tau$ are inverse isomorphisms of $C(K*L)$ and $D(K*L).$
Moreover, the action $t_\sharp$ on $C(K*L)$, defined by $t_K$,$t_L$, corresponds to the action $t_\sharp(\sigma \otimes \tau) = t_{K\sharp}(\sigma) \otimes t_{L\sharp}(\tau)$ on $D(K*L).$
\end{proposition}

Note that $\emptyset*\sigma=\sigma$ and $\tau*\emptyset = \tau$ for any simplices $\sigma, \tau,$ in particular, $\emptyset * \emptyset$ is the -1-dimensional cell of $K*L.$ The above proposition is proved by straightforward calculations using the standard formula for the boundary of $\sigma * \tau.$ Let $d_\sigma$, $d_\tau$ denote the dimensions of these simplices. If $d_\sigma,d_\tau \geq 0$ then the standard boundary formula implies,

$$\partial (\sigma *\tau) = \sum_{i=0}^{d_\sigma} (-1)^i \widehat{\sigma_i}*\tau + (-1)^{d_\sigma+1} \sum_{j=0}^{d_\tau} (-1)^j \sigma*\widehat{\tau_j},$$
where $\widehat{\sigma_i}$ is the oriented simplex obtained from removing the $i$-th vertex of $\sigma$, and $\widehat{\tau_j}$ is defined analogously. If one of $d_\sigma$ and $d_\tau$ is negative or both are, then the above formula still works if we make the convention that $\sigma*0=0*\tau=0$ for any simplices $\sigma$, $\tau.$

\subsection{Resolutions in the Join}
Denote by $\phi_i, i\in {0,1, \ldots}$ an arbitrary resolution of $1 \in C^0_d(K)$ and by $\phi'_i, i \in \{0,1, \ldots\}$ an arbitrary resolution of $1 \in C^0_d(L).$ Also we write $\Phi_i, i\in \{0,1,\ldots\}$ for the resolution of $1 \in D^0_d(K*L)$ which we are going to compute. Recall that $\emptyset \in K$ is the $-1$-dimensional simplex and the generator of $C_{-1}(K)$ and that for any vertex $v$, $\partial v = \emptyset.$ We also denote by $\emptyset \in C^{-1}(K, \Zring)$ the dual cochain to $\emptyset$ which assigns 1 to $\emptyset.$\\

\paragraph{Notation} In the rest of this paper we denote the action on both $K$ and $L$ by $t$. Moreover, for simplicity we write $t$ for $t^\sharp$.\\

\paragraph{}
Note that for all $i,j \in \Zring$ there are isomorphisms 
$$\mu: \text{Hom}(C_i(K), \Zring) \otimes \text{Hom}(C_j(L), \Zring) \cong  \text{Hom}(C_i(K) \otimes C_j(K),\Zring)  $$
given by the homomorphism $\mu(\eta \otimes \eta')(\sigma \otimes \tau) = \eta(\sigma)\eta'(\tau)$, see for instance \cite[Section 5.5]{Spa89}. Therefore there are isomorphisms
$$\mu:  \sum_{j} \text{Hom}(C_{n-1-j}(K), \Zring) \otimes \text{Hom}(C_j(L), \Zring) \cong  \text{Hom}(D_n(K*L),\Zring).$$
We use these isomorphisms implicitly.

\begin{lemma}\label{l:joinreszero}
We can set $\Phi_0= \emptyset \otimes \phi'_0 + \phi_0 \otimes \emptyset.$
\end{lemma}

\begin{proof}
We need to check that $s \Phi_0=1$. We have $s \Phi_0 = \emptyset \otimes \phi'_0 + \phi_0 \otimes \emptyset + \emptyset \otimes t \phi'_0 + t\phi_0 \otimes \emptyset + \cdots + \emptyset \otimes t^{p-1} \phi'_0 + t^{p-1}\phi_0 \otimes \emptyset.$ Note that a ``vertex" in $D_0(K*L)$ is of the form $\emptyset \otimes v$ or $w\otimes \emptyset.$ Since $s\phi_0=1$ and $s\phi'_0=1$, on any vertex exactly one of the terms evaluates to 1 and the rest are zeros. It follows that $s \Phi_0=1.$ 
\end{proof}

In order to compute $\Phi_i,$ $i>0,$ we simply need to successively apply the definition of a resolution. We present here the result of such a calculation and we defer the tedious derivation of these formulas to Section~\ref{s:rescalc}.

\begin{proposition}\label{p:joinres}
Let $\phi_i, \phi'_i, i\in \{0,1,2, \ldots\}$ be resolutions of $1$ for $(K,t)$ and $(L,t)$ respectively. A resolution
$\Phi_j$ of $1$ for $K*L$ is defined as follows. For any $j\in \{0, 1, 2, \ldots\},$

$$ \Phi_{j} = \Phi^{j}_{-1} + \Phi^{j}_0 + \Phi^{j}_{1}+ \cdots + \Phi^{j}_{j},$$
where,

\begin{enumerate}
    \item If $j=2m$ is even then
    
    \begin{align*}
    \Phi^{2m}_{-1} &= \emptyset \otimes \phi'_{2m}\\
    \Phi^{2m}_{2l} &= \phi_{2l} \otimes t^{-l} \phi'_{2m-2l-1}\\
    \Phi^{2m}_{2l+1} &= \phi_{2l+1} \otimes t^{-(l+1)} \phi'_{2m-(2l+1)-1}.
    \end{align*}
    
    \item If $j=2m+1$ is odd, set $s_q = 1+t+ \ldots+t^q$, $q \in \{0,1,\ldots,p\}$, then, 

    \begin{align*}
    \Phi^{2m+1}_{-1} &= \emptyset \otimes \phi'_{2m+1}\\
    \Phi^{2m+1}_{2l} &= \sum_{q=0}^{p-2}{\left((s-s_q) \phi_{2l}\right) \otimes t^{-l+q}\phi'_{2m+1-2l-1}}\\
    \Phi^{2m+1}_{2l+1} &= \phi_{2l+1} \otimes t^{-(l+1)} \phi'_{(2m+1)-(2l+1)-1}.
    \end{align*}    
    
\end{enumerate}

\end{proposition}


\paragraph{Remark}
Of course, in Proposition \ref{p:joinres} the cochains $\phi_i$, $\phi'_i$, $\Phi_j$ will be 0 if their index is larger than the dimension of the complex. Also note that $\Phi^j_k$ involves $\phi_k$, $\phi'_{k'}$ such that $k+k'=j-1.$ \\

\paragraph{Remark}
We have utilized augmented complexes and chain complexes in order to have the structure to prove Proposition~\ref{p:joinres}. After we have computed a chain complex and a resolution for a join, then we can replace the $-1$-dimensional group with 0 and work with the non-augmented chain complex.  

\paragraph{Remark}
Note that for $p=2$ the formula in the proposition remains unchanged other than for the term $\Phi_{2l}^{2m+1}$ which simplifies to

$$ \Phi_{2l}^{2m+1}= t\phi_{2l} \otimes t^{-l} \phi'_{2m+1-2l-1}.$$

\section{The Smith Index of the Join}\label{s:joinsmith}
In this section, we prove our main theorems that determine the Smith index of the join of two $\Zring_p$-complexes in terms of the indices of the factors. Recall that for a $\Zring_p$-complex $K$ we denote the $i$-dimensional Smith class by $A^i(K)=A^i(K,t)$, and when computed using modulo $p$ operations we denote the resulting class with $A^i_p(K)$. Also recall that by the top Smith class we always mean the Smith class (with the coefficients group understood from the context) of the largest dimension which is non-zero. The dimension of the top Smith class is one less than the Smith index. 

We first prove a useful lemma.
\begin{lemma}[Shortening of the Smith resolution]\label{l:shortening}
Let $K$ be a $\Zring_p$-complex with Smith index $I(K)$. There is a resolution $\phi_i$, $i=1, \ldots$ of $1 \in C^0(K)$ such that $\phi_j=0$ for $j \geq I(K)$.
\end{lemma}

\begin{proof}
Assume $I(K)$ is even. By assumption we have $A^{I(K)}(K)=0$. Therefore, $s\phi_{I(K)} = \delta s \psi$ for some cochain $\psi \in C^{I(K)-1}(K).$ Let $\phi_i$ be an arbitrary resolution of $1\in C^0(K)$. We modify the resolution $\phi_i$ into $\overline{\phi}_i$ by setting $\overline{\phi}_i = \phi_i$ for $i<I(K)-1$, $\overline{\phi}_{I(K)-1} = \phi_{I(K)-1}- s\psi$, and $\overline{\phi}_i =0 $ for $i> I(K)-1.$ Then $\delta \overline{\phi}_{I(K)-2} = d \phi_{I(K)-1} = d (\phi_{I(K)-1}-s\psi)$ and the new sequence $\overline{\phi}_i$ is a resolution of $1$ for $K.$ We do totally analogously if $I(K)$ is odd.

\end{proof}

\begin{theorem}\label{t:upperbound}
Let $K$ and $L$ be $\Zring_p$-complexes. Let $I(K)$ and $I(L)$ be the Smith indices of $K$ and $L$. Then 
$A^{I(K)+I(L)}(K*L)=0$ and therefore 

\begin{equation}\label{eq:upperbound}
    I(K*L) \leq I(K)+I(L).
\end{equation}

\end{theorem}

\begin{proof}
Let $\phi_i$ and $\phi'_i$ be resolutions of $1$ for $K$ and $L$ respectively. 
Applying Lemma~\ref{l:shortening} we assume that the resolutions are shortened, that is,
$\phi_{I(K)}=0$ and $\phi'_{I(L)}=0$.

Using Proposition \ref{p:joinres} we can construct a resolution of 1 for $K*L$ using the new resolutions, denoted $\overline{\Phi}_j.$ Consider now a term $\overline{\Phi}_j$ where $j\geq I(K)+I(L).$ Recall that $\overline{\Phi}_j = \sum_{k=-1}^{j} \overline{\Phi}^j_k.$ Now  $\overline{\Phi}^j_k$ will be zero if $k \geq I(K)$ or $j-1-k\geq I(L).$ One of these happens if $j-1 \geq I(K)+I(L)-1$, which is our assumption.

\end{proof}

\begin{definition}
If equality happens in (\ref{eq:upperbound}) then we say that the Smith index is \textit{stable} under join of $K$ and $L$. 
\end{definition}


\subsection{Conditions for the Stability of the Smith Index }
We next prove the stability of the Smith index under certain conditions. We first consider the case where one of the join factors $K$ or $L$ is a set of $p$ points. This case is known for $p=2,$ see \cite[ (5.1)]{CoFl60}, and is called the \textit{stability under suspension}. Our notion of stability expands the stability under suspension to joins with any complex. In the case of joins with 0-dimensional complexes, we prove more than stability so that later we can use it to simplify our arguments.

\begin{definition}
Let $K$ be a $\Zring_p$-complex and $n_K$ be the dimension of $K$'s top Smith class. A sequence of numbers $p, p^{m_1}, \ldots, p^{m_{n_K}}$ is called the \textit{sequence of moduli} for $K$ if for each $1 \leq j \leq n_K,$ $A^j(K)\neq 0 \mod p^{m_j}$ and if $m_j>1$, $A^j(K)= 0 \mod p^{m_j-1}.$
\end{definition}

\begin{lemma}\label{l:zerodim}
Let $\Sigma_p= \{ a_1, \ldots, a_p\}$ be a discrete set of $p$ points and set $t(a_i)=a_{i+1}$, where $a_{p+1}=a_1$. Let $K$ be any $\Zring_p$-complex. Let $p,p^{m_1}, \ldots,p^{m_{n_K}}$ be the sequence of moduli for $K$. Then the sequence of moduli for $\Sigma_p*K$ is $p,p,p^{m_1}, \ldots, p^{m_{n_K}}.$ In particular, 
$$I(\Sigma_p * K)= I(K)+1, \; I_p(\Sigma_p*K) = I_p(K)+1.$$
\end{lemma}



\begin{proof}
Let $\phi_i$ be a resolution of 1 for $\Sigma_p$ and $\phi'_i$ be a resolution of 1 for $K.$ We use the shortened resolution  of Lemma~\ref{l:shortening} for $K$. Thus $\phi'_{I(K)}=0.$ Let $\Phi_0, \Phi_1, \ldots$ be a resolution of 1 for $\Sigma_p*K$. Let $n \geq 1$ be an integer. We consider two cases: $n$ odd and even. 

Assume $n$ is odd. We first show that $A^{n-1}(K)=0$ modulo $p^m$, for some $m,$ implies $A^{n}(\Sigma_p*K)=0$ modulo $p^m.$ If $n=1$ then the condition implies that $A^0(K)=0$ and $K$ is empty. We therefore assume $n \geq 3.$ The assumption implies that there are cochains $\psi$ and $x$ such that $s\phi'_{n-1} = \delta s \psi +p^m sx.$
We change the term $\phi'_{n-2}$ of the resolution to $\bar{\phi}'_{n-2}=\phi'_{n-2}-s \psi.$ Then still $\delta \phi'_{n-3}=d \bar{\phi}'_{n-2}$ and $\delta \bar{\phi}'_{n-2}= \delta \phi'_{n-2}-\delta s \psi =s\phi'_{n-1}-\delta sp= p^m s x$. Thus we can set $\bar{\phi}'_{n-1}=p^m x.$ And we obtain $\delta \bar{\phi}'_{n-1}= p^m \delta x = d \bar{\phi}'_{n}.$
Using Proposition \ref{p:joinres}, we have
\begin{align*}
\Phi_{n} &= \emptyset \otimes \bar{\phi}'_{n} + \sum_{q=0}^{p-2}((s-s_q)\phi_0) \otimes t^q \bar{\phi}'_{n-1}\\
&= \emptyset \otimes \bar{\phi}'_{n} \mod p^m.\\
\end{align*}
However, $d \Phi_{n} = d(\emptyset \otimes \bar{\phi}'_{n}) = \emptyset \otimes d \bar{\phi}'_{n}=0$, modulo $p^m.$ Therefore, $A^{n}(\Sigma_p*K)=0$ modulo $p^m.$

To finish the proof in this case, we need to provide a certificate that the $s$-cohomology class of $d\Phi_{n}$ is non-trivial modulo $p^{m_{n-1}}$. We can take $\phi_0$ such that $\phi_0(a_1)=1$ and $\phi_0(a_i)=0$, $i\neq 1$. Let $c'$ be the chain given as a certificate by Lemma~\ref{l:certd} for $A^{n-1}(K)$ with modulus $p^{m_{n-1}}$ with respect to a fundamental domain of the action of $\Zring_p$ which contains the support of $\phi'_{n-1}$. Here we assume that $\phi'_{n-1}$ is chosen such that its support lies in a single fundamental domain and $s\phi'_{n-1}=\delta \phi'_{n-2}$.

Consider the chain $da_1\otimes sc'$. Note that $da_1 \otimes sc' = d(a_1 \otimes sc')$ and $da_1 \otimes sc'$ is an $s$-chain. We claim that $x=a_1 \otimes sc'$ is a certificate for $d\Phi_{n}$ modulo $p^{m_{n-1}}.$ We compute that $\partial(dx)= \partial (da_1 \otimes sc')= \partial (da_1) \otimes sc' - da_1 \otimes \partial (sc') = 0-da_1\otimes s\partial (c')= 0 \mod p^{m_{n-1}}.$ In addition, the only summand of $\Phi_n$ which can be non-zero on $da_1\otimes sc'$ is $\Phi^n_0$ and we have
$\Phi_{n}(da_1 \otimes sc')=\sum_{q=0}^{p-2} (s-s_q)\phi_0(da_1)\cdot \phi'_{n-1}(sc').$
Now one checks that $\sum_{q=0}^{p-2}(s-s_q)d = s-p\cdot id$. Therefore $\sum_{q=0}^{p-2} (s-s_q)\phi_0(da_1) = \phi_0(sa_1-pa_1)=1-p$. Thus $\Phi_{n}(da_1 \otimes sc') = (1-p)\phi'(sc') \neq 0 \mod p^{m_{n-1}}.$ The latter number is non-trivial modulo $p^{m_{n-1}}$ since $\phi'_{n-1}$ is zero outside the fundamental domain and $c'$ is a linear combination of oriented simplices of the fundamental domain hence $\phi(sc')=\phi(c').$

Next we consider the case where $n$ is even and greater than zero.

We first show that if $A^{n-1}(K)=0$ modulo $p^m$, for some $m,$ then $A^{n}(\Sigma_p*K)=0$ modulo $p^m.$ The assumption implies that there are cochains $\psi$ and $x$ such that $d\phi'_{n-1} = \delta d \psi +p^m dx.$
We change the term $\phi'_{n-2}$ of the resolution to $\bar{\phi}'_{n-2}=\phi'_{n-2}-d \psi.$ Then still $\delta \phi'_{n-3}=s \bar{\phi}'_{n-2}$ (if $n\geq 3$) and $\delta \bar{\phi}'_{n-2}= \delta \phi'_{n-2}-\delta d \psi = p^m d x$. Thus we can set $\bar{\phi}'_{n-1}=p^m x.$ And we obtain $\delta \bar{\phi}'_{n-1}= p^m \delta x = s \bar{\phi}'_{n}.$
Using Proposition \ref{p:joinres}, we have
\begin{align*}
\Phi_{n} &= \emptyset \otimes \bar{\phi}'_{n} + \phi_0 \otimes \bar{\phi}'_{n-1}\\
&= \emptyset \otimes \bar{\phi}'_{n} \mod p^m.\\
\end{align*}
However, $s \Phi_{n} = s(\emptyset \otimes \bar{\phi}'_{n}) = \emptyset \otimes s \bar{\phi}'_{n}=0$, modulo $p^m.$ Therefore, $A^{n}(\Sigma_p*K)=0$ modulo $p^m.$

For the rest of the proof in this case, we need to provide a certificate that the $d$-cohomology class of $s\Phi_{n}$ is not zero modulo $p^{m_{n-1}}.$ Since $n$ is even Proposition~\ref{p:joinres} tells us that $\Phi_{n}$ contains the summand
$\Phi^{n}_0=\phi_0 \otimes \phi'_{n-1}.$ The cochain $\Phi^{n}_0$ is zero outside a fundamental domain containing $\{a_1*\sigma\}$ where $\sigma$ ranges over all $n-1$-simplices of $K$. Let $c'$ be a certificate for $A^{n-1}(K)$ with $n=p^{m_{n-1}}.$ Consider the chain $a_1 \otimes dc'.$ We have $s \partial (a_1 \otimes dc')= s(\partial (a_1) \otimes dc')-s(a_1 \otimes \partial (dc')) = \emptyset \otimes sdc' - s(a_1 \otimes \partial (dc'))=0 \mod p^{m_{n-1}}.$ On the other hand, $(\phi_0 \otimes \phi'_{n-1})(a_1\otimes dc')= \phi_0(a_1)\phi'_{n-1}(dc')=\phi'_{n-1}(dc') \neq 0 \mod p^{m_{n-1}}$.

\end{proof}


\begin{definition}
Let $K$ be $\Zring_p$-complex. Let $h$ be a $d$-cohomology class and $c$ be a certificate for $h$ with modulus $n$. Let $q>0$ be an integer. We say $c$ is \textit{boundary-equivariant mod $q$} if $d\partial c=0 \mod q$, or equivalently, $t\partial c= \partial c \mod q$.

\end{definition}

Note that if $p=2$ every certificate is boundary-equivariant modulo 2 since $d=s$ modulo 2. 
The condition implies that the boundary of the certificate is an equivariant chain modulo $q$. In other words, the boundary of the certificate is equal modulo $q$ to some chain in the image of the transfer homomorphism, if the quotient complex is defined. As an example, an even-dimensional Smith class in an antipodal sphere has a dual which is boundary-equivariant, that is, whose boundary is equivariant (in integers). This is simply a half-sphere.

Recall that when computing modulo $p$, each non-trivial cocycle has a dual cycle, which is a $p$-cycle on which the cocycle is non-trivial modulo $p$. We call a chain $c$ \textit{dual} to a $d$-cocycle $s\phi$ if $sc$ is dual to $\gamma(s\phi)$. Analogously, we define duality for $s$-cycles and $s$-chains. Note that the certificates given by Lemmas~\ref{l:certd}~and~\ref{l:certs} for $n=p$ are dual to the given homology class. When working modulo $p$ we use certificate and dual interchangeably.

\begin{theorem}\label{t:stabil}
Let $K$ and $L$ be $\Zring_p$-complexes and $n_K, n_L$ be positive numbers.

\begin{enumerate}

    \item If $p=2$ and $A^{n_K}(K) \neq 0 \mod 2$ and $A^{n_L}(L) \neq 0 \mod 2$, then $A^{n_K+n_L+1}(K*L) \neq 0 \mod 2.$
    
    
    \item If $n_K>0$ is even, $A^{n_K}(K) \neq 0 \mod p$ and $A^{n_L}(L) \neq 0 \mod p^m$, for some minimal $m\geq 1$, then $A^{n_K+n_L+1}(K*L) \neq 0 \mod p^m$ if $A^{n_K}(K)$ has a dual which is boundary-equivariant modulo $p^m$. 
    
\end{enumerate}

\end{theorem}

\begin{proof}

The case where $K$ or $L$ is 0-dimensional is covered by Lemma~\ref{l:zerodim}. Therefore, we assume $K$ and $L$ have positive dimensions. Let $A=A^{n_K}(K)$ and $A'=A^{n_L}(L).$ As usual, assume $\phi_i$ and $\phi'_i$ are resolutions of 1 for $K$ and $L$ respectively. 


Since for $p=2$ every certificate is $2$-boundary-equivariant, statement (1) is a special case of (2). However, we give a simpler proof for this case.
First we prove statement (1), when both $n_K$ and $n_L$ are even and thus $n_K+n_L+1$ is odd. 
We can assume that the supports of $\phi_{n_K}$ and $\phi'_{n_L}$ are over a fundamental domain. The cochain $s \phi_{n_K}$ is a representative for the $d$-cohomology class $A$. Since $A_p$ is non-trivial by the assumption, we can take the $n$ given by Lemma~\ref{l:certd} for $A$ to be $p$. Therefore there is a chain $c \in C_{n_K}(K)$ such that $sc$ is a cycle modulo $p$ and $\phi_{n_K}(c) \neq 0$ modulo~$p$. Similarly a certificate $c'$ exists for $A'_p=[s \phi'_{n_L}]_d$ such that $sc'$ is a cycle modulo $p$ but $\phi_{n_L}(c')\neq 0 \mod p.$

In the rest of the proof we write $\phi$ for $\phi_{n_K}$ and ${\phi}'$ for ${\phi'}_{n_L}$ and we set $l=(n_K-1)/2.$ We need to present a certificate for $A^{n_K+n_L+1}(K*L)$ in the sense of Lemma~\ref{l:certs}. 

Assume $p=2$, and note that in this case $s=d \mod 2$. We claim that the chain $c\otimes sc'$ is a certificate. We first compute,

\begin{align*}
    \partial (d (c\otimes sc') ) &= d \left( \partial(c) \otimes sc' - c \otimes \partial (sc') \right)\\
    &= s\partial(c) \otimes sc' \mod 2\\
    &=0 \mod 2.
\end{align*}
Thus our candidate for a certificate is a $d$-cocycle. In the rest of the proof we do not assume
$p=2$ for generality. We do this to emphasize that we need $p=2$ to satisfy the condition on the boundary of our certificate, and that this boundary condition is what goes wrong for $p>2$.

Using Proposition \ref{p:joinres} we know that the only term of $\Phi_{n_K+n_L+1}$ which is possibly non-zero on $c \otimes sc'$ is $\Phi^{n_K+n_L+1}_{n_K}$. Recall that a term of the form $\eta^i \otimes \mu^j$, by the definition of the tensor product of cochains, is only non-zero on the subspace of the chains of the form $x \otimes y$ where $x$ is an $i$-chain and $y$ is a $j$-chain. 
Since $n_K$ is even and $n_K+n_L+1$ is odd, the middle case of part (2) of Proposition \ref{p:joinres} gives

\begin{equation}
   \Phi^{n_K+n_L+1}_{n_K} = \sum_{q=0}^{p-2}{\left((s-s_q) \phi_{n_K}\right) \otimes t^{-l+q}\phi'_{n_L}}.
\end{equation}

Using Lemma~\ref{l:claim1} below with $x=\phi_{n_K}$ and $y=t^{-l}\phi'_{n_L}$, and noting that $s\phi'_{n_L}=st^{-l}\phi'_{n_L}$, we obtain 

\begin{equation}\label{eq:dphi}
   d\Phi^{n_K+n_L+1}_{n_K} = s\phi \otimes t^{-l}\phi' - \phi \otimes s \phi'. 
\end{equation} 

We now need to show that $d\Phi^{n_K+n_L+1}_{n_K} (c\otimes sc')\neq 0 \mod p$.
The second term on the right assigns a multiple of $p$ to $c \otimes sc'$ and the first term assigns the number $\phi(c)\phi(c')$ which is non-zero modulo $p$. Therefore the theorem is proved in this case.

Now one can reduce the other cases of statement (1) to the above using Lemma~\ref{l:zerodim}. For instance, assume $n_K$ and $n_L$ are odd. Then we replace $K$ with $K'=S^0*L$ and $L$ with $L'=S^0*L$ so that we are in the case above. Lemma~\ref{l:zerodim} implies that $A^{n_L+1}(L')\neq 0$ modulo $p$ and $A^{n_K+1}(K') \neq 0$ modulo $p$. The above argument proves that $A^{n_K+n_L+3}(K'*L')\neq 0$ modulo $p$. On the other hand, by two applications of Lemma~\ref{l:zerodim}, $A^{n_K+n_L+3}(S^0*S^0*K*L)\neq 0$ modulo $p$ implies $A^{n_K+n_L+1}(K*L)\neq 0$ modulo $p$ as was required.

We now prove statement (2). We assume there exists a certificate $c$ (with $n=p$) for $A^{n_K}(K)$ such that $\partial c$ is boundary-equivariant modulo $p^m$, i.e.,  $\partial c=s x \mod p^m$ for some chain $x$, and that $A^{n_L}(L)$ has a certificate modulo $p^m$ and $m$ is minimal.

First assume $n_K$ and $n_L$ are both even. Again, we define the certificate for $A^{n_K+n_L+1}(K*L)$ to be $c \otimes sc'$. We compute

\begin{align*}
    \partial (d (c\otimes sc') ) &= d \left( \partial(c) \otimes sc' - c \otimes \partial (sc') \right)\\
    &= d \left( s x \otimes sc' - c \otimes \partial(sc') \right) \mod p^m\\
    &=dsx \otimes sc' \mod p^m\\
    &=0 \mod p^m.
\end{align*}

To calculate the value of $\Phi^{n_K+n_L+1}_{n_K}$ on $d(c\otimes sc')$ we again use formula~(\ref{eq:dphi}). Since $m$ is minimal and Smith classes are $p$-torsion, $\phi'(sc')=p^{m-1}y$ for some $y\neq 0 \mod p$ and the second term on the right gives $\phi(c)\cdot s\phi'(sc)$ and is zero modulo $p^m$. The first term is $\phi(sc)\cdot\phi'(sc')$ and is non-zero modulo $p^m$, as was required.

\end{proof}

\begin{lemma}\label{l:claim1}

Let $x \in \text{Hom}(C_i(K), \Zring), y \in \text{Hom}(C_j(L), \Zring)$ be arbitrary cochains. Then,
\begin{align*}
    d \left( \sum_{k=0}^{p-2} (s-s_k)x \otimes t^k y \right) = sx\otimes y - x \otimes sy.
\end{align*}
\end{lemma}

\begin{proof}
\begin{align*}
    &(1-t) \left( \sum_{k=0}^{p-2} (s-s_k)x \otimes t^k y \right)\\ 
    &=  \sum_{k=0}^{p-2} (s-s_k)x\otimes t^k y - \sum_{k=0}^{p-2} (s-s_k)tx \otimes t^{k+1}y\\
    &= - \sum_{k=0}^{p-2} s_k x\otimes t^k y + \sum_{k=0}^{p-2} s_k tx \otimes t^{k+1}y +\sum_{k=0}^{p-2} sx \otimes (t^k - t^{k+1}y)\\
    &= -s_0 x \otimes t^0 y + \sum_{k=1}^{p-2} (-s_k x + s_{k-1}t x) \otimes t^k y + s_{p-2}tx \otimes t^{p-1} y \\ &+ sx \otimes \sum_{k=0}^{p-2}  (t^k - t^{k+1})y \\
    &= -x \otimes y - \sum_{k=1}^{p-2}  x \otimes t^k y + s_{p-2}t x \otimes t^{p-1} y + sx \otimes (1- t^{p-1}) y \\
    &=  -x \otimes y - \sum_{k=1}^{p-2}  x \otimes t^k y + (s_{p-2}t - s)x \otimes t^{p-1} y + sx \otimes y \\
    &= -x \otimes sy + sx \otimes y.
\end{align*}

\end{proof}

\subsection{Examples}\label{examples}
\begin{enumerate}[a)]
    \item\label{ex:a} Here is an example of a $\Zring_2$-complex whose top Smith class is non-trivial as an integer class but trivial modulo $2$. Take a 1-sphere $S^1$ with the antipodal action. Glue two 2-cells to the sphere just as to make a 2-sphere but instead map the boundaries by degree 2 maps. Extend the antipodal action of $S^1$ by mapping one disk to the other. Now a triangulation of this space is the $\Zring_2$-complex with the property we are looking for, see Figure~\ref{fig:examplea}. Observe that the result is no longer a 2-manifold. One can also obtain such complexes as the deleted products or deleted joins of Melikhov complexes defined in Section~\ref{s:embed}.
    
    \begin{figure}[ht]
        \centering
        \includegraphics[scale=0.4]{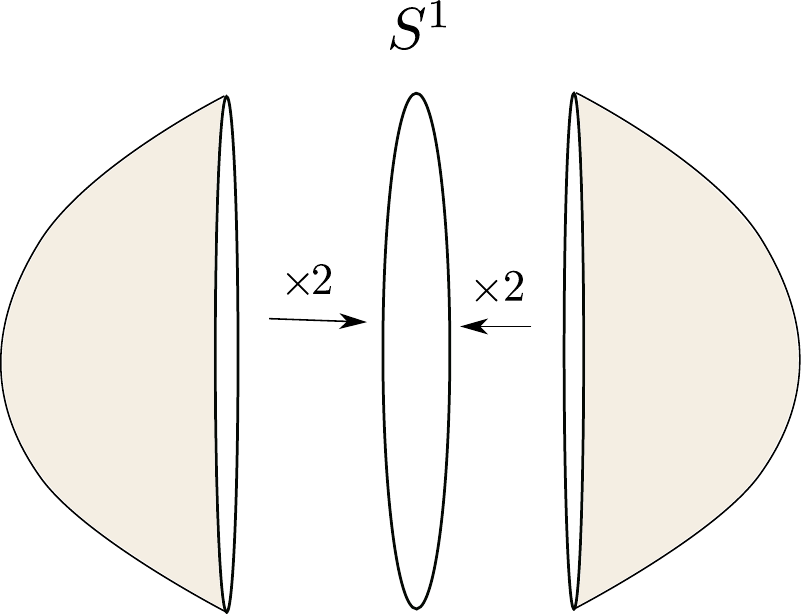}
        \caption{The complex of Example~\ref{ex:a}}
        \label{fig:examplea}
    \end{figure}
    
    \item\label{ex:b} We present a $\Zring_2$-complex whose top Smith class is non-zero modulo 2, moreover, it does not satisfy the condition of statement (2) of Theorem \ref{t:stabil}. Take a 2-sphere with the antipodal action consisting of two 2-cells $D_1$ and $D_2$. We remove interiors of two small disks $E_{11}$ and $E_{12}$ from $D_1$ and we remove the interiors of the images of these under the antipodal action from $D_2$, namely $E_{21}$ and $E_{22}$, respectively. We then glue a cylinder $I \times S^1$ to $\partial E_{11}$ and $\partial E_{22}$ such that one boundary of the cylinder is mapped to $\partial E_{11}$ by a degree 1 map, and the other boundary is mapped to $\partial E_{22}$ by a map of degree $-1$. Analogously, we glue in a cylinder to $\partial E_{12}$ and $\partial E_{22}$. We extend the antipodal map to the resulting space by mapping one cylinder to the other. See Figure~\ref{fig:exampleb}. Let $P$ be a triangulation of this space which defines a $\Zring_2$-complex. The 2-dimensional Smith class of this $\Zring_2$-complex $P$ does not have a dual which is boundary-equivariant modulo $2^2$. Our computer calculations show that $P$ indeed has non-vanishing modulo~2 Smith class of dimension~2. Moreover, we have calculated the 5-dimensional Smith class of the join of $P$ with the complex of Example~\ref{ex:a} and it is indeed zero. Therefore, this example shows that in general the Smith index of the join is not determined by the indices of the factors alone. 
    
    \begin{figure}[ht]
        \centering
        \includegraphics[scale=0.4]{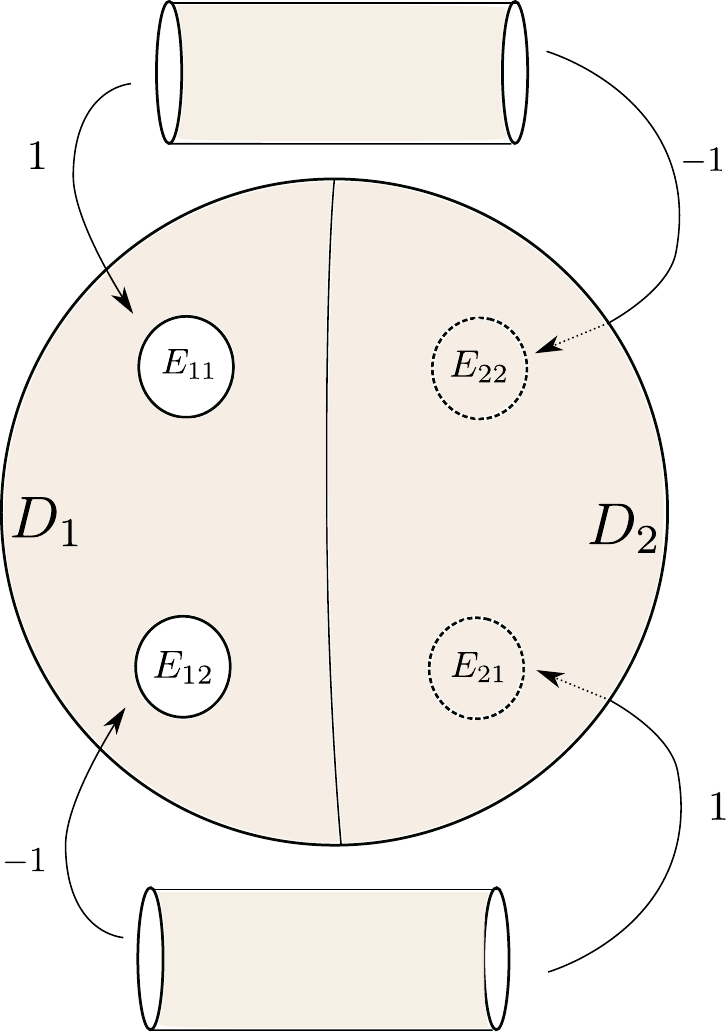}
        \caption{The complex of Example~\ref{ex:b}}
        \label{fig:exampleb}
    \end{figure}    
    
\end{enumerate}

\subsection{Instability of the Smith Index}
We now prove that if for both $K$ and $L$ the modulo $p$ reduction of the top Smith class is trivial then the Smith index is not stable.

We start by proving a lemma which plays an essential role in the argument.

\begin{lemma}\label{l:eqresterm}
Let $K$ be a simplicial $\Zring_p$-complex such that $I_p(K)=I(K)-1$. In other words, the top non-zero Smith class of $K$ is trivial modulo $p$. If $I(K)-1$ is even then there is a resolution $\phi_i, i=1,\cdots$ such that $\phi_{I(K)-1}=s\psi$ for some cochain $\psi$, and $\phi_j=0$, for $j>I(K)-1$.
\end{lemma}

\begin{proof}
Let $A:=A(K)^{I(K)-1}$. Assume $I(K)-1$ is even. Since the classes $A$ is zero modulo $p$ but non-zero as an integer class we can write $A=pB$ for non-zero classes $B.$ Let $s\eta$ be a $d$-cocycle representing $B$. We modify an arbitrary resolution $\phi_i$ into $\tilde{\phi}_i$ and $\tilde{\phi}'_i$ such that $\tilde{\phi}_{I(K)-1}=s\eta$ as follows. For simplicity we do not write the subscripts $I(K)-1$ and we write $-k$ for $I(K)-1-k$ and $+k$ for $I(K)-1+k.$ 

For $j < I(K)-2$ we do not change the terms $\phi_j$. Since $s\phi-sp\eta = \delta(sb)$ for some cochain $b,$ we set $\tilde{\phi}_{-1}= \phi_{-1} - sb.$ Then $\delta (\tilde{\phi}_{-1}) = s \phi - \delta(sb)=ps\eta=ss\eta=s\tilde{\phi}.$ And still $\delta(\tilde{\phi}_{-2})=d \tilde{\phi}_{-1}.$  Since $\delta (s\eta)=0$ the higher terms of the resolution are zeros. 
\end{proof}

\begin{theorem}\label{t:instab}
Let $K$ and $L$ be simplicial $\Zring_p$-complexes such that $I_p(K)=I(K)-1$ and $I_p(L)=I(L)-1.$ Then $I(K*L) \leq I(K)+I(L)-1$. Moreover, the Smith class $A^{I(K)+I(L)-2}(K*L)$ is zero modulo $p$.
\end{theorem}

\begin{proof}
We first consider the case where $I(K)-1$ and $I(L)-1$ are both even. Let as usual $\phi_i$ be a resolution of 1 for $K$ and $\phi'_i$ be a resolution of 1 for $L.$ The assumption is that $A=A^{I(K)-1}(K)$ and $A'=A^{I(L)-1}(L)$ both vanish modulo $p.$ 

Apply Lemma~\ref{l:eqresterm} to obtain resolutions $\tilde{\phi}_i$ and $\tilde{\phi'}_i$ such that $\tilde{\phi}_{I(K)-1}=s\eta$ and $\tilde{\phi}'_{I(L)-1} = s \eta'$ for some cochains $\eta$ and $\eta'$. 

Let $2l = I(K)-1$. Using Proposition \ref{p:joinres}, we know that the only summand of $\Phi_{I(K)+I(L)-1}$ which could possibly be non-zero is $\Phi^{I(K)+I(L)-1}_{2l}$. We compute,

\begin{align*}
    \Phi^{I(K)+I(L)-1}_{2l}&= \sum_{q=0}^{p-2}(s-s_q)\tilde{\phi}\otimes t^{-l+q}  \tilde{\phi}' \\
    &= \sum_{q=0}^{p-2} (s-s_q)s\eta\otimes t^{-l+q} s \eta'\\
    &= \left(\sum_{q=0}^{p-2} (s-s_q)s \eta \right) \otimes s\eta'\\
    &= \frac{p(p-1)}{2} s\eta \otimes s \eta'.
\end{align*}
On the other hand the Smith class of dimension $I(K)+I(L)-1$ of $K*L$ is the $s$-cohomology class of
$d \Phi^{I(K)+I(L)-1}_{2l}.$ However,

\begin{align*}
    d\Phi^{I(K)+I(L)-1}_{2l} &= \frac{p(p-1)}{2} d (s\eta \otimes s\eta')\\
    &= \frac{p(p-1)}{2} ds\eta \otimes s \eta' = 0.
\end{align*}
It follows that $I(K*L) \leq I(K)+I(L)-1.$


Next we prove the last assertion of the theorem. We could have defined a resolution, Smith classes, etc. using chain complexes with $\Zring_p$ coefficients instead
of the integer coefficients. Let $A_{(p)}^j$ be the Smith classes defined using $\Zring_p$ coefficients. Then the $\Zring_p$ versions of Proposition~\ref{p:joinres} and Theorem~\ref{t:upperbound} would be true if $p$ is prime. These then would imply that $A^{I(K)+I(L)-2}_{(p)}(K*L)=0.$ To prove the assertion, we still need to show that $A^j_{(p)}$ classes are equal to the reductions modulo $p$ of the integer Smith classes. This follows since the resolution with coefficients $\Zring_p$ can be taken to be the reduction of the integer resolution. To see this latter statement, observe that it is true for the 0-th term of the resolution, and it follows for higher terms since the reduction modulo $p$ homomorphisms commute with $\delta$ and $t.$ This finishes the proof for the case $I(K)-1$ and $I(L)-1$ both even. 

Entirely analogously to the proof of Theorem~\ref{t:stabil}, Lemma~\ref{l:zerodim} now allows us to reduce other cases to the above case.
\end{proof}

Since the stability result of Theorem~\ref{t:stabil} involves the Smith classes with non-zero reduction modulo $p$ the following is very useful. 

\begin{lemma}\label{l:topmod}
Let $K$ be a $\Zring_p$-complex. Only the top Smith class of $K$ can possibly be both non-trivial as an integer class and trivial modulo $p.$
\end{lemma}

\begin{proof}

The simplest argument to prove the statement uses Lemma~\ref{l:eqresterm}. Let $m$ be the smallest integer such that $A^m(K)$ is zero when computed modulo $p$. Using Lemma~\ref{l:zerodim} we can assume $m$ is even. Then we can replace the $m$-th term of the resolution with a cocycle. Implying that higher terms of the resolution can be set to zeros.

A more intuitive argument uses the cup product structure of Smith classes and the fact that these classes are $p$-torsion elements. This structure is described in \cite{Sha57} for $p=2.$ If $A^{n_K-1}$ is zero modulo $2$ then we can write $A^{n_K-1}=2A'.$ Then $A^{n_K}=A^1 \smallsmile A^{n_K-1} = A^1 \smallsmile 2 A' = 2A^1 \smallsmile A'=0,$ which is a contradiction.

\end{proof}

\begin{corollary}\label{c:topm}
For any $\Zring_p$-complex $K,$ the sequence of moduli is of the form $p,p, \ldots, p, p^m,$ for some $m \geq 1.$
\end{corollary}

We can now give a proof of Theorem~\ref{t:smith}.
\begin{proof}[Proof of Theorem~\ref{t:smith}]
Statement (1) follows from Theorem~\ref{t:instab}. Statement (2) follows from (1) of Theorem~\ref{t:stabil} by setting $n_K=I_2(K)-1$ and $n_L=I_2(L)-1$ and Theorem~\ref{t:upperbound}.

To prove statement (3) assume $A^{(I_p(K)-1)}(K) \neq 0 \mod p^m$ for some minimal $m$. 
Note that $I(K*L)\geq I(K)+I(L)-1$ since in Theorem~\ref{t:stabil}, with Lemma~\ref{l:topmod} in mind, we can take $n_K=I(K)-2$ and $n_L=I(L)-1$, in which case the theorem implies that the Smith class of dimension $I(K)+I(L)-2 \neq 0$ modulo 2. If $A^{I_2(L)-1}(L)$ has a dual which is boundary-equivariant mod $2^m$ then the index is stable. This happens for instance when $L$ is an antipodal 2-sphere. Example b) gives a $\Zring_2$-complex $L$ whose top Smith class does not have a dual with equivariant boundary and our computations show that the Smith index is not stable under join with this complex.

By ``extra conditions" in statement (4) we actually mean the conditions in (2) of Theorem~\ref{t:stabil} when $n_K=I(K)-1$ and $n_L=I(L)-1$.

\end{proof}

\section{Embeddability of Joins of Complexes}\label{s:embed}
In this section, we prove our theorems on the embeddability of the join $M*N$ of two simplicial complexes
$M$ and $N.$ We only consider embeddability in the double dimension. That is, we answer the questions of the form: Is there an embedding of an $n$-dimensional simplicial complex into the $2n$-dimensional Euclidean space? Our aim is to derive conditions for the embeddability of the join complex $M*N$ into double dimension based on the properties of $M$ and $N$ and their respective embeddability or non-embeddability into double dimension.

For $n=1,$ the problem of embeddability into double-dimension is deciding planarity of graphs, which, needless to say, is very well-known. Flores~\cite{Flo32} and van Kampen~\cite{vKam33} were the first to construct examples of $n$-dimensional simplicial complexes that do not embed into $2n$-dimensional space for any $n>1$. The method used by van Kampen was developed and completed by Shapiro~\cite{Sha57} and Wu~\cite{Wu74}, giving rise to the van Kampen obstruction class and the theory of embedding classes. Note that we refer only to the embedding class of twice the dimension as \textit{the} van Kampen obstruction class, and by embedding classes we refer to this or lower dimensional cohomological obstructions to embeddability. The van Kampen class is usually defined as a cohomology class of the quotient of the deleted product of the complex whose vanishing is a necessary and sufficient condition for a complex of dimension $n>2$ to be embeddable into $2n$-space. For $n=2,$ this obstruction is known to be not complete~\cite{Freetal94} and no algorithm is known for deciding the embeddability of a 2-complex into 4-space via a PL map.  See \cite{Matetal10} for very interesting hardness results in this direction.

The van Kampen obstruction can be defined in several ways. For the purposes of this paper it is best to define it as a Smith class as follows. Let $M$ be a $d_M$-dimensional simplicial complex. Consider the \textit{deleted product} $M^{\times 2}_\Delta$ of $M$ which is, by definition, the subcomplex of $M \times M$ consisted of products of vertex-disjoint simplices of $M.$ In other words, remove from $M \times M$ every cell $\sigma \times \sigma'$ if $\sigma \cap \sigma' \neq \emptyset.$ 
There is an action of $\Zring_2$ on $M^{\times 2}_\Delta$ whose non-trivial map $t$ sends $(x,y) \in |M^{\times 2}_\Delta|$ onto $(y,x)$. This action simply permutes the cells and is determined by the action on vertices, namely, $t(v,w)=(w,v)$ where $v$ and $w$ are distinct vertices of $M.$ Since the deleted product consists of cells which are not simplices we can either take a simplicial sub-division of it or we can temporarily allow ourselves to consider polyhedral $\Zring_p$-complexes. The definition of Smith classes is given analogously, in a general complex. Since we do not apply our results to the deleted product, we allow ourselves this temporary extension of our definition of a $\Zring_p$-complex. Therefore, we have a $\Zring_2$-complex $(M^{\times 2}_\Delta,t)$. We define the van Kampen obstruction class of $M$ to be the $d$-cohomology class $A^{2d_M}(M^{\times 2}_\Delta).$ This definition is justified by the theorem of Section V.5 of \cite{Wu74} which asserts the Smith class and the van Kampen obstruction with the traditional definition coincide. Intuitively, since $2d_M$ is even, they are both images of the same cohomology class of the space $\Rspace P^{\infty}$ under the induced homomorphism of any $\Zring_2$-equivariant map $f:M^{\times 2}_\Delta \rightarrow S^{\infty},$ where $S^\infty$ is the standard infinite dimensional sphere with the antipodal action. Note that, here we have also changed the traditional definition by defining the van Kampen obstruction to be a $d$-cohomology class of the deleted product rather than an ordinary class of the quotient space. Lemma~\ref{l:dquotientiso} says that these two cohomology groups are isomorphic in even dimensions.

We summarise some known results on the properties of van Kampen obstruction class and conditions for the embeddability in double dimension.
\begin{proposition}
Let $M$ be a simplicial complex of dimension $d_M$.

\begin{enumerate}
    \item If $d_M \geq 2$ then there exists a $\Zring_2$-equivariant map $f: M^{\times 2}_\Delta \rightarrow S^{2d_M-1}$ if and only if the van Kampen obstruction vanishes, i.e., $A^{2d_M}(M^{\times 2}_\Delta)=0.$
    
    \item If $d_M \geq 3$ then there exists a PL embedding $g:M \rightarrow \Rspace^{2{d_M}}$ if and only if there exists a $\Zring_2$-equivariant mapping $f: M^{\times 2}_\Delta \rightarrow S^{2d_M-1}.$
\end{enumerate}
\end{proposition}

Statement (1) follows using the ``finger moves" argument for $d_M\geq 3$ and in general using Lemma~\ref{l:obst} below. Statement (2) expresses the sufficiency of the van Kampen obstruction class or is a special case of the Haefliger-Weber criteria  for embedding in the meta-stable range \cite{Web67}.

One can define the obstructions to embedding in the more amenable deleted join rather than deleted product. However, the classical sufficiency theorems work with the deleted product and the properties of the deleted join have to be translated to the deleted product for the purposes of demonstrating an embedding.

The \textit{deleted join} of a simplicial complex $M$, denoted $M^{*2}_\Delta$, is the subcomplex of the simplicial complex $M*M$ consisted of the joins of vertex-disjoint simplices. The action of $\Zring_2=\{ 1,t \}$ on $M^{*2}_\Delta$ is the simplicial action defined by $t(v*\emptyset) = \emptyset*v$ and $t(\emptyset*v) = v*\emptyset$. If a complex $M$ embeds into $\Rspace^{2d_M}$ then there is an equivariant embedding of $M^{*2}_\Delta$ into the ``deleted join" of $\Rspace^{2d_M}$. The deleted join can also be defined for general spaces, see \cite{Mat08} for details. For our purposes, it is enough to recall that if $M$ embeds into $R^{2d_M}$ then there is an equivariant map $M^{*2}_\Delta \rightarrow S^{2d_M}$ \cite[Lemma~5.5.4]{Mat08}, and therefore $I(M^{*2}_\Delta) \leq 2d_M+1.$ Therefore, the Smith class $A^{2d_M+1}(M^{*2}_\Delta)$ is an obstruction for embeddability into double dimension.

We gather here some important lemmas for future reference. For the proof of the following lemma see \cite[ part (e) of Lemma 4]{PaSk20}.
\begin{lemma}\label{l:joinpro}
If $d_M \geq 2$, then there exists a $\Zring_2$-equivariant map $f_\times : M^{\times2}_\Delta \rightarrow S^{2d_M-1}$ if and only if there exists a $\Zring_2$-equivariant map $f_*: M^{*2}_\Delta \rightarrow S^{2d_M}.$
\end{lemma}

\begin{lemma}\label{l:obst}
Let $K$ be a $\Zring_2$-complex of dimension $n>3$. Then the class $A^{n}(K)$ is a complete obstruction for existence of a $\Zring_2$-equivariant map $K \rightarrow S^{n-1}.$
\end{lemma}

\begin{proof}
The necessity of $A^{n}(K)=0$ for existence of the map is easy. 
First assume $n=2m$. The proof of the sufficiency in this case can be found in \cite{Mel09} proof of Theorem 3.2. We expound on that proof. One can define the map on $(2m-1)$-skeleton arbitrarily. This will give a map from the $(2m-1)$-skeleton of the quotient complex $K/t$ to $\Rspace P^\infty.$ If we can extend this latter map to all of $K/t$ and lift it to an equivariant map we are done. The lift is possible if $2m-1>1$ or $m>1.$ The complete obstruction to the extension problem is the preimage of the single element in the ordinary cohomology class of dimension $2m$ of $\Rspace P^{\infty}$ which is the same as the isomorphic image of the Smith class of dimension $2m$ in the homology of the quotient space. This follows for instance by the theory of obstruction classes developed in \cite{Ste51}. In our case the bundle of coefficients is trivial, see \cite[37.6]{Ste51}. 

If $n=2m+1$ then we first suspend to obtain a complex $\Sigma(K)$ of even dimension. Lemma~\ref{l:zerodim} says that $A^n(K)\neq 0$ if and only if $A^{n+1}(\Sigma(K))\neq 0$.
Now by the equivariant suspension theorem~\cite{CoFl60}[Theorem 2.5], if $n \leq 2n-5$ then the existence of an $\Zring_2$-equivariant map from $\Sigma(K)$ into $S^{n}$ implies the existence of a $\Zring_2$-equivariant map from $K$ into $S^{n-1}$. Thus if $A^n(K)$ is trivial there is an equivariant map into $S^{n-1}$. If there is an equivariant map into $S^{n-1}$ then $A^n(K)$ has to be trivial.

\end{proof}

\subsection{Main Theorems on the Embeddability of Join of Complexes}\label{ss:embed}
We now present the application of our results of Section \ref{s:joinsmith} to the problem of embeddability of joins of complexes in double dimension.

\begin{theorem}\label{t:joinembed}
Let $M$ and $N$ be simplicial complexes of dimensions $d_M$ and $d_N,$ respectively, such that $d_M+d_N\geq 2$. Assume $I(M^{*2}_\Delta)=2d_M+2$ and $I(N^{*2}_\Delta)=2d_N+2$. Thus $M$ and $N$ are not embeddable into $\Rspace^{2d_M}$ and $\Rspace^{2d_N}$ respectively. 

\begin{enumerate}
    \item If the top Smith classes of $M^{*2}_\Delta$ and of $N^{*2}_\Delta$ are non-zero modulo 2 then $M*N$ does not embed into the Euclidean space of dimension $2(d_M+d_N+1)$.
    
    \item If one of the top Smith classes of $M^{*2}_\Delta$ or of $N^{*2}_\Delta$, say $M^{*2}_\Delta$, is non-zero modulo 2 and has a dual which is boundary-equivariant modulo $2^m$, and the top Smith class of $N^{*2}_\Delta$ is non-zero modulo $2^m$ with $m$ minimal,  then $M*N$ does not embed into the Euclidean space of dimension $2(d_M+d_N+1)$.
    
    \item If the top Smith classes of $M^{*2}_\Delta$ and of $N^{*2}_\Delta$ are both zero modulo 2 then $M*N$ embeds into the Euclidean space of dimension $2(d_M+d_N+1)$.
\end{enumerate}
\end{theorem}

\begin{proof}
Assume first that the reduction modulo 2 of the classes $A^{2d_M+1}(M^{*2}_\Delta)$ and $A^{2d_N+1}(N^{*2}_\Delta)$ are zeros.
Theorem \ref{t:instab} then says that $A^{2d_M+2d_N+3}( M^{*2}_\Delta * N^{*2}_\Delta)=0.$ There is a well-known $\Zring_2$-homeomorphism $M^{*2}_\Delta * N^{*2}_\Delta \cong (M*N)^{*2}_\Delta,$ see for instance \cite[Lemma 5.5.2]{Mat08}. Therefore, $A^{2d_M+2d_N+3}((M*N)^{*2}_\Delta)=0.$ Lemma~\ref{l:obst} then implies that there is an equivariant map $(M*N)^{*2}_\Delta \rightarrow S^{2d_M+2_N+2}.$  Next, by Lemma~\ref{l:joinpro} there is an equivariant map $(M*N)^{\times 2}_\Delta \rightarrow S^{2d_M+2d_N+1}.$ Therefore $M*N$ embeds in the $(2d_M+2d_N+2)$-dimensional Euclidean space since $d_M+d_N+1 \geq 3$. This proves statement (3).

Assume, on the other hand, that the assumptions of statements (1) or (2) are satisfied. Then by Theorem \ref{t:stabil}, $I((M*N)^{*2}_\Delta)= I(M^{*2}_{\Delta}*N^{*2}_{\Delta})=I(M^{*2}_{\Delta})+ I(N^{*2}_{\Delta}) = 2d_M+2d_N+4.$  It follows that $A^{2d_M+2d_N+3}((M*N)^{*2}_\Delta)\neq 0$ and therefore $M*N$ does not embed into $(2(d_M+d_N+1))$-dimensional Euclidean space.

\end{proof}

We need also show that there exist complexes $M$ and $N$ that satisfy the conditions of Theorem~\ref{t:joinembed}. Simplicial complexes whose deleted join has vanishing modulo 2 Smith class of dimension 5 are given by Melikhov \cite{Mel09}, see \ref{ss:melikhov} below. We can now give the proof of Theorem~\ref{t:existence}.

\begin{proof}[Proof of Theorem~\ref{t:existence}]
As shown in \cite{Mel09} for any $d$ there is a complex $M$ of dimension $d$ whose van Kampen obstruction is zero modulo 2 but non-zero as an integer class, see the following subsection for a description of these complexes. We need to show that the $(2d+1)$-dimensional Smith class of the deleted join $M^{*2}_\Delta$ of these complexes is non-zero but zero modulo 2. 


Note that $M^{*2}_\Delta$ is equivariantly homeomorphic to $(CM)^{\times 2}_\Delta,$ where $CM$ is the cone over $M,$ see \cite[Exercise 4 to 5.5]{Mat08}. Also, as shown in \cite{Sko02}, there exists an equivariant surjective map $p: (CM)^{\times 2}_\Delta \rightarrow \Sigma M^{\times 2}_\Delta$ where $\Sigma$ denotes suspension. The composition map $M^{*2}_\Delta \cong_{\Zring_2} (CM)^{\times 2}_\Delta \xrightarrow{p} \Sigma M^{\times 2}_\Delta$, induces an isomorphism on homologies of dimension greater than $d_M+1$, since $p$ contracts two sub-complexes of dimensions $d_M\geq 2$ into points. The claim now follows since the suspension changes the sequence of moduli as in Lemma~\ref{l:zerodim}.


\end{proof}

We gather the results on the conditions for the embeddability of $K*L$ into the double dimension.

\begin{theorem}\label{t:mainembed}
Let $M$ and $N$ be simplicial complexes of dimensions $d_M$ and $d_N$ respectively. 

\begin{enumerate}

    
    \item The complex $M*N$ embeds into the $2(d_M+d_N+1)$-dimensional Euclidean space if the van Kampen obstruction of one of $M$ or $N$ vanishes or both have vanishing obstructions modulo~2. 
    
    \item The complex $M*N$ does not embed into the $2(d_M+d_N+1)$-dimensional Euclidean space if both obstructions are non-zero and one of the following is satisfied.
    \begin{enumerate}
        \item $d_M \leq 1$ or $d_N \leq 1$.
        \item The van Kampen obstructions of $M$ and $N$ are non-zeros modulo 2.
        \item One of the van Kampen obstruction classes is non-zero modulo 2 and the corresponding Smith class in the deleted join has a dual which is boundary-equivariant mod $2^m$, and the other van Kampen obstruction class is non-zero modulo $2^m$ with $m$ minimal. 
    \end{enumerate}

\end{enumerate}

\end{theorem}

\begin{proof}
Observe that the van Kampen obstruction of a complex cannot be non-trivial yet zero modulo 2 unless the complex is at least 2-dimensional. This is because for a graph vanishing of the van Kampen obstruction modulo 2 implies vanishing integer obstruction class. This follows since $K_5$ and $K_{3,3}$ have non-zero modulo 2 obstructions.

Assume $d_M+d_N=0$. Then $M$ and $N$ are discrete point sets and
their obstruction classes are non-zero if and only if they contain at least three points. If one of $K$ or $L$ is a two-point set then $M*N$ clearly embeds in $\Rspace^2.$ This proves statement (1) in this case.

Assume $d_M+d_N=1$, then say $M$ is a discrete point set and $N$ is a graph. If $M$ contains two points, then $M*N$ is suspension of a graph that embeds into $\Rspace^4$ since any graph embeds into $\Rspace^3.$ If the graph $N$ is planar it is easily seen that $M*N$ embeds into $\Rspace^4.$ This proves (1) in this case. 


Assume $d_M+d_N\geq 2$ and that $M$ has vanishing obstruction. If $M$ is zero-dimensional then as before $M*N$ embeds in double dimension. 
If $M$ is 1-dimensional, i.e. a graph, then the vanishing of the van Kampen obstruction is equivalent to $M$ being planar. On the other hand, $N$ embeds into $\Rspace^{2d_N+1}.$ Therefore, the join $M*N$ can be embedded into the join $\Rspace^2*\Rspace^{2d_N+1} = \Rspace^{2d_M+2d_N+2}.$ Note that this argument could not be used if $M$ were 2-dimensional, since then vanishing of the van Kampen obstruction does not imply existence of an embedding into $\Rspace^4.$

If $M$ is of dimension $> 1$, then $M^{\times 2}_\Delta$ is of dimension $2d_M > 2$. Since the van Kampen obstruction of $M$ vanishes, by Lemma~\ref{l:obst}, there is an equivariant map $M^{\times 2}_\Delta \rightarrow S^{2d_M-1}.$ Then by Lemma~\ref{l:joinreszero} there is an equivariant map $M^{* 2}_\Delta \rightarrow S^{2d_M}.$ It follows that $I(M^{* 2}_\Delta) \leq 2d_M+1.$
Theorem \ref{t:upperbound} then implies that $I((M*N)^{*2}_\Delta) \leq I(M^{* 2}_\Delta)+I(N^{*2}_\Delta) \leq 2d_M+1+2d_N+2.$ From Lemma~\ref{l:obst} it follows that there is an equivariant map $(M*N)^{* 2}_\Delta \rightarrow S^{2d_M+2d_N+2}$. Lemma~\ref{l:joinpro} then gives an equivariant map $(M*N)^{\times 2}_\Delta \rightarrow S^{2d_M+2d_N+1}$. Finally $M*N$ embeds into $\Rspace^{2d_M+2d_N+2}$ since $d_M+d_N+1 \geq 3$ and its van Kamen obstruction vanishes. 
This proves statement (1) in the case $d_M+d_N\geq 2$. 

Assume one of $M$ or $N$, say $M$, is of dimension less than 2. Since the obstruction of $M^{*2}_\Delta$ is non-trivial, it is non-trivial modulo 2 and the join $M^{*2}_\Delta$ contains an antipodal sphere. From either (repeated application of) Lemma~\ref{l:zerodim} or Theorem~\ref{t:stabil}.2 it follows that the Smith index is stable under the join of the deleted joins $M^{*2}_\Delta$ and $N^{*2}_\Delta$ and the join is non-embeddable into the double dimension as in the proof of Theorem~\ref{t:joinembed}. This proves (2.a). 

The remaining statements follow from Theorem \ref{t:joinembed}. 


\end{proof}

Since the top Smith class of an even-dimensional sphere with the antipodal action has a dual whose boundary is equivariant, there are many simplicial complexes whose van Kampen obstructions have this property. For instance, if the deleted product contains an antipodal sphere as an equivariant subset. In fact, we do not know if there exists a simplicial complex whose obstruction class is non-trivial modulo 2 but does not have a dual with equivariant boundary, see Section~\ref{s:discussion}.

\subsection{Melikhov Complexes}\label{ss:melikhov}
In \cite[Example 3.6]{Mel09}, Melikhov first answered the question of whether the modulo 2 van Kampen obstruction is a complete obstruction whenever the integer class is, by constructing complexes whose modulo 2 obstruction vanishes but not the integer obstruction. These are quite simple and here we elucidate this construction. For any $n\geq 2,$ let $\Delta^{(n)}_{2n+2}$ be the $n$-skeleton of the $(2n+2)$-dimensional simplex. For example, when $n=2,$ $\Delta^{(2)}_6$ is the collection of all triangles and edges that one can build from seven vertices (together with their faces). Van Kampen \cite{vKam33} has shown that $\Delta^{(n)}_{2n+2}$ is not embeddable in the Euclidean space of double dimension. It is also not difficult to construct a simplex-wise linear map of $\Delta^{(n)}_{2n+2}$ into $\Rspace^{2n}$ with a single transverse intersection between (images of) vertex-disjoint simplices.

\begin{figure}
    \centering
    \includegraphics[scale=0.6]{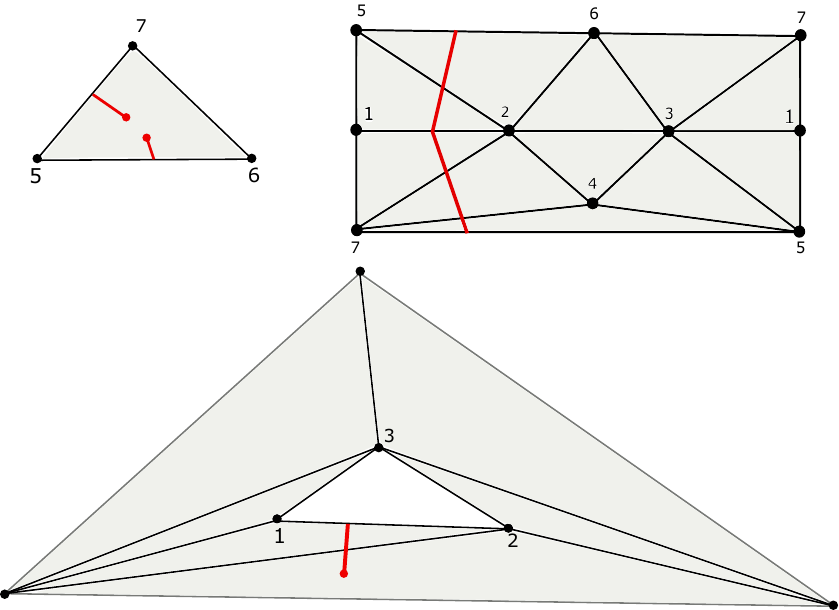}
    \caption{Constructing a complex with vanishing mod 2 obstruction, one triangle in the 2-skeleton of the 6-simplex is replaced by the result of identifying vertices with the same name in this figure}
    \label{fig:mel}
\end{figure}

Now the Melikhov complexes are built by taking a single $n$-simplex of $\Delta^{(n)}_{2n+2}$, removing a small open disk from it, and gluing back the disk by a map of degree 2. In the case $d=2$,  the complex that replaces a triangle of $\Delta^{(2)}_6$ is depicted in Figure~\ref{fig:mel}. This complex is obtained by gluing a capped Mobius band to the boundary of the removed disk, along the central circle of the band. A visual proof that the modulo 2 obstruction of the resulting complex vanishes is also given in the figure. Without loss of generality, we can assume that the triangle $T$ which we are replacing, contains the single intersection point. This point is depicted with a red dot in the figure. By pushing the (image of the) triangle that intersects $T$ we can push the intersection point to the cap of the Mobius band, such that it creates a red trace on the complex, which looks like a fork, as in the figure. At the end of this move, there are two same-sign intersection points between two vertex disjoint simplices and no other intersections. It follows that the van Kampen obstruction class has a cocycle which is zero modulo 2. Recall that the van Kampen obstruction class is traditionally defined using the intersection homomorphism, in the quotient of the deleted product,  which assigns the algebraic intersection number to pairs of disjoint top-dimensional simplices. Here this representative of the obstruction class gives 2 on a single cell and is zero elsewhere. The proof that the obstruction is not zero as an integer can be found in \cite{Mel09}. 

Observe that if we replace the Mobius band in this construction by a tower of Mobius bands consisted of $k$ bands, then we can make the obstruction vanish modulo $2^{k}$ but not modulo $2^{k+1}.$ 

\section{Calculation of a Resolution for the Join Complex}\label{s:rescalc}
In this section we prove Proposition \ref{p:joinres}. We need some lemmas that we now state and prove.
Throughout this section, $\phi_i$ is a resolution of $1$ for the $\Zring_p$-complex $K$ and $\phi'_i$ is a resolution of $1$ for the $\Zring_p$-complex $L$. Recall that we denote the action of $t \in \Zring_p = \{1,t, \ldots, t^{p-1}\}$ on $K$ and $L$ and also their join by $t$. No confusion should arise. Moreover for simplicity of notation as usual we denote the homomorphism $t^\sharp$ by $t$.

\begin{lemma}\label{l:claim2}
Let $x \in \text{Hom}(C_i(K), \Zring), y \in \text{Hom}(C_j(L), \Zring)$ be arbitrary cochains. Then,
\begin{align*}
    x \otimes sy + \sum_{k=0}^{p-2}\left( (s-s_k)dx\right) \otimes t^k y = s\left(x \otimes t^{p-1}y\right)
\end{align*}

\end{lemma}

\begin{proof}
Since $(s-s_k)d = t^{k+1}-1$ we obtain,
\begin{align*}
    &x \otimes sy + \sum_{k=0}^{p-2}\left( (s-s_k)dx \right) \otimes t^k y \\
    &= x \otimes sy + \sum_{k=0}^{p-2} \left( (t^{k+1} -1)x \right)\otimes t^k y  \\
    &= x \otimes t^{p-1} y + \sum_{k=0}^{p-2} t^{k+1} x \otimes t^k y \\
    &= x \otimes t^{p-1} y + t^{k+1} \sum_{k=0}^{p-2} x \otimes t^{p-1} y = s\left(x \otimes t^{p-1}y \right).
\end{align*}
\end{proof}

\begin{lemma}\label{l:claim3}

Let $x \in \text{Hom}(C_i(K), \Zring), y \in \text{Hom}(C_j(L), \Zring)$ be arbitrary cochains. Then,
\begin{align*}
    sx \otimes y - \sum_{k=0}^{p-2}\left( (s-s_k)x \right) \otimes t^k d y = s \left(x \otimes y \right).
\end{align*}
\end{lemma}

\begin{proof}
We simplify the second term on the left,

\begin{align*}
&-\sum_{k=0}^{p-2}( (s-s_k)x )\otimes t^k d y\\
&= \sum_{k=0}^{p-2} \left( -sx \otimes t^k(1-t)y + s_kx \otimes t^k (1-t)y \right)\\
&= \sum_{k=0}^{p-2} \left( -t^k(sx \otimes y) + t^k (sx \otimes ty) + t^k\left(s_k t^{-k} x \otimes y\right) - t^k\left(s_kt^{-k}x \otimes ty \right)   \right)\\
&= \sum_{k=0}^{p-2} \left(  -t^k(sx \otimes y) +t^{k+1}(sx \otimes y) + t^k\left(s_k t^{-k} x \otimes y \right) - t^{k+1}\left(s_{k}t^{-(k+1)}x \otimes y \right) \right)\\
&= - sx \otimes y + t^{p-1} (sx \otimes y)\\&+  \sum_{k=0}^{p-2} \left( t^k \left(s_k t^{-k} x \otimes y\right) - t^{k+1}\left(s_{k+1}t^{-(k+1)}x \otimes y\right) + t^{k+1}(x \otimes y) \right)\\
&= - sx \otimes y +  t^{p-1}(sx \otimes y) + x\otimes y -t^{p-1}\left(s_{p-1}t^{-(p-1)}x \otimes y\right) \\
&+ \sum_{k=0}^{p-2} t^{k+1}(x \otimes y)\\&= - sx \otimes y + s(x \otimes y). 
\end{align*}

\end{proof}

\begin{lemma}\label{l:claim4}

Let $x \in \text{Hom}(C_i(K), \Zring), y \in \text{Hom}(C_j(L), \Zring)$ be arbitrary cochains. Then,
\begin{align*}
    dx \otimes y + x \otimes t^{p-1} d y = (1-t)\left(x \otimes t^{p-1} y\right).
\end{align*}
\end{lemma}

\begin{proof}
We compute,
\begin{align*}
    dx \otimes y + x \otimes t^{p-1} d y &= \left( (1-t)x \right) \otimes y + x \otimes\left( t^{p-1}y - t^p y  \right)\\
    &= x \otimes y - tx\otimes y + x\otimes t^{p-1}y - x \otimes y\\
    &= -t ( x\otimes t^{p-1}y ) + x\otimes t^{p-1}y\\
    &= (1-t)\left(x \otimes t^{p-1} y\right).
\end{align*}

\end{proof}

Now using the above lemmas we derive the formulas of Proposition \ref{p:joinres} for the resolution of the join complex. We prove that the given formulas are correct for the 0-dimensional term of the resolution and then we prove that given the formula for the term of dimension $n$ we obtain the formula for the term of dimension $n+1.$

Proposition \ref{p:joinres} gives $\Phi_0 = \emptyset \otimes \phi'_0 + \phi_0 \otimes \emptyset.$ Note that for any $x,y$, $t(\emptyset \otimes y) = \emptyset \otimes (ty)$ and $t(x \otimes \emptyset) = (ty) \otimes \emptyset.$ From these relations it follows that $s(\Phi_0)= \emptyset \otimes (s \phi'_0) + (s\phi_0) \otimes \emptyset = \emptyset \otimes 1 + 1 \otimes \emptyset.$ This latter cochain assigns 1 to any vertex of $K*L.$

Now let $n=2m$ be even. We compute,

\begin{align*}
    \delta\left(\Phi^{2m}_{-1}\right)&= \delta \left(\emptyset \otimes \phi'_{2m}\right)\\
    &= \delta\left(\emptyset\right) \otimes \phi'_{2m}+\left(-1\right)^{0} \emptyset \otimes \delta\left(\phi'_{2m}\right)\\
    &= s\phi_0 \otimes \phi'_{2m} + \emptyset \otimes d\phi'_{2m+1},\\
    \delta\left(\Phi^{2m}_{2l}\right)&= \delta\left(\phi_{2l} \otimes t^{-1} \phi'_{2m-2l-1}\right)\\
    &= \delta\left(\phi_{2l}\right) \otimes t^{-l} \phi'_{2m-2l-1} + \left(-1\right)^{2l+1} \phi_{2l} \otimes t^{-l}\delta\left(\phi'_{2m-2l-1}\right)\\
    &=d\phi_{2l+1}\otimes t^{-l}\phi'_{2m-2l-1} - \phi_{2l}\otimes t^{-l}s\phi'_{2m-2l}\\
    &=d\phi_{2l+1}\otimes t^{-l}\phi'_{2m-2l-1} - \phi_{2l}\otimes s\phi'_{2m-2l},\\
    \delta\left(\Phi^{2m}_{2l+1}\right) &= \delta\left(\phi_{2l+1} \otimes t^{-\left(l+1\right)} \phi'_{2m-2l-2}\right)\\
    &= \delta\left(\phi_{2l+1}\right) \otimes t^{-\left(l+1\right)} \phi'_{2m-2l-2} + \left(-1\right)^{2l+2} \phi_{2l+1} \otimes t^{-\left(l+1\right)}\delta\left(\phi'_{2m-2l-2}\right)\\
    &= s\phi_{2l+2}\otimes t^{-\left(l+1\right)}\phi'_{2m-2l-2} + \phi_{2l+1} \otimes t^{-(l+1)}d\phi'_{2m-2l-1}.
\end{align*}

There are 6 terms (including the preceding signs) in the right hand sides of the above and we number them from left to right and top to bottom. We now collect the terms with the same dimension in the first factor, for various $l$, to form the terms of $d\Phi_{2m+1}$. We start at dimension $k=-1$ where we have only the term 
\begin{align*}
  \emptyset \otimes d \phi'_{2m+1} &= d(\emptyset \otimes \phi'_{2m+1}).  
\end{align*}
For the dimension $k=0$ we have the term number 1 above together with the term number 4 when $l=0.$ Then using Lemma~\ref{l:claim1} we obtain
\begin{align*}
    s\phi_0 \otimes \phi'_{2m} - \phi_0 \otimes s \phi'_{2m} = d \left(\sum_{q=0}^{p-2} (s-s_q)\phi_0 \otimes t^{q}\phi'_{2m}\right).
\end{align*}
 If the dimension $k>0$ is even then we put together again the term number 4 when $l=k/2$ and the term number~5 when $l=(k-2)/2.$ Again using Lemma~\ref{l:claim1} we obtain
 \begin{align*}
    s\phi_k \otimes t^{-((k-2)/2+1)} \phi'_{2m-k} -& \phi_k \otimes s t^{k/2} \phi'_{2m-k}\\
    &= s\phi_k \otimes t^{-k/2} \phi'_{2m-k} - \phi_k \otimes s t^{k/2} \phi'_{2m-k}\\
    &=d \left(\sum_{q=0}^{p-2} (s-s_q)\phi_k \otimes t^{q-(k/2)}\phi'_{2m-k}\right).
\end{align*}

If $k$ is odd then we have the term number 3 with $l=(k-1)/2$ and term number 6 again with $l=(k-1)/2$. Using Lemma~\ref{l:claim4} we obtain
 \begin{align*}
    d \phi_{k} \otimes t^{-(k-1)/2} \phi'_{2m-k} &+ \phi_k \otimes t^{-((k-1)/2)-1} d\phi'_{2m-k}\\
    &=d \phi_k \otimes t^{-(k-1)/2} \phi'_{2m-k} + \phi_k \otimes d t^{p-1} t^{-(k-1)/2} \phi'_{2m-k}\\
    &=d\left(\phi_k \otimes t^{p-1} t^{-(k-1)/2} \phi'_{2m-k}\right)\\
    &=d\left(\phi_k \otimes t^{-(k-1)/2-1} \phi'_{2m-k}\right).
\end{align*}

By comparing the terms we have computed for each $k$ with the formula claimed by Proposition \ref{p:joinres} we see that if $\Phi_{2m}$ is as given by the proposition then the formula claimed for $\Phi_{2m+1}$ is correct.

We next consider the case where $n=2m+1$ is odd. We compute
\begin{align*}
    \delta\left(\Phi^{2m+1}_{-1}\right) &= \delta \left(\emptyset \otimes \phi'_{2m+1}\right)\\
    &= \delta\left(\emptyset\right) \otimes \phi'_{2m+1}+\left(-1\right)^{0} \emptyset \otimes \delta\left(\phi'_{2m+1}\right)\\
    &= s\phi_0 \otimes \phi'_{2m+1}+ \emptyset \otimes s\phi'_{2m+2},\\   
    \delta\left(\Phi^{2m+1}_{2l}\right) &= \delta\left(\sum_{q=0}^{p-2}\left(\left(s-s_q\right)\phi_{2l}\right) \otimes t^{q-l}\phi'_{2m-2l}\right)\\
    &= \sum_{q=0}^{p-2} \left[ \left(s-s_q\right) \delta \left(\phi_{2l}\right) \otimes t^{q-l}\phi'_{2m-2l}+\left(-1\right)^{2l+1}\left(s-s_q\right)\phi_{2l} \otimes t^{q-l}\delta\left(\phi'_{2m-2l}\right) \right]\\
    &= \sum_{q=0}^{p-2} \left[ \left(s-s_q\right) d\left(\phi_{2l+1}\right) \otimes t^{q-l}\phi'_{2m-2l}-\left(s-s_q\right)\phi_{2l} \otimes t^{q-l}d\phi'_{2m-2l+1} \right]\\
    &= \sum_{q=0}^{p-2} \left( -s_qd\phi_{2l+1} \otimes t^{q-l}\phi'_{2m-2l}\right) - \sum_{q=0}^{p-2}\left(s-s_q\right)\phi_{2l} \otimes t^{q-l}d \phi'_{2m-2l+1}.
\end{align*}
We simplify the first term of the right hand side above,
\begin{align*}
     &\sum_{q=0}^{p-2} \left( -s_qd\phi_{2l+1} \otimes t^{q-l}\phi'_{2m-2l}\right)\\ 
     &=  \sum_{q=0}^{p-2} \left( \left(t^{q+1}-1\right)\phi_{2l+1} \otimes t^{q-l}\phi'_{2m-2l}\right)\\
      &=\sum_{q=0}^{p-2} \left[ t^{q+1}\phi_{2l+1} \otimes t^{q-l}\phi'_{2m-2l} - \phi_{2l+1}\otimes t^{q-l}\phi'_{2m-2l} \right]\\
      &= \sum_{q=0}^{p-2}t^{q+1}\left(\phi_{2l+1}\otimes t^{-l-1}\phi'_{2m-2l}\right) - \phi_{2l+1}\otimes \left(\sum_{q=0}^{p-2}t^{q-l}\phi'_{2m-2l}\right)\\
      &= \left(s-1\right)\left(\phi_{2l+1} \otimes t^{-l-1}\phi'_{2m-2l}\right) - \phi_{2l+1}\otimes \left(s-t^{p-1}\right)t^{-l}\phi'_{2m-2l}\\
      &= s\left(\phi_{2l+1}\otimes t^{-l-1}\phi'_{2m-2l}\right)-\phi_{2l+1}\otimes t^{-l-1}\phi'_{2m-2l} - \phi_{2l+1}\otimes s\phi'_{2m-2l}\\&+\phi_{2l+1}\otimes t^{-l-1}\phi'_{2m-2l}\\
      &= s\left(\phi_{2l+1}\otimes t^{-l-1}\phi'_{2m-2l}\right)-\phi_{2l+1}\otimes s\phi'_{2m-2l}.
\end{align*}
We thus obtain,
\begin{align*}
    \delta\left(\Phi^{2m+1}_{2l}\right) &= s\left(\phi_{2l+1}\otimes t^{-l-1}\phi'_{2m-2l}\right)-\phi_{2l+1}\otimes s\phi'_{2m-2l}\\ &- \sum_{q=0}^{p-2}\left(s-s_q\right)\phi_{2l} \otimes t^{q-l}d \phi'_{2m-2l+1}.
\end{align*}

We compute the coboundary in the last case,
    \begin{align*}
        \delta\left(\Phi^{2m+1}_{2l+1}\right) &= \delta\left(\phi_{2l+1}\otimes t^{-l-1}\phi'_{2m-2l-1}\right)\\
        &=\delta\left(\phi_{2l+1}\right)\otimes t^{-l-1}\phi'_{2m-2l-1} +\left(-1\right)^{2l+2}\phi_{2l+1}\otimes t^{-l-1}\delta\left(\phi'_{2m-2l-1}\right)\\
        &=s\phi_{2l+2}\otimes t^{-l-1} \phi'_{2m-2l-1} + \phi_{2l+1}\otimes s\phi'_{2m-2l}.
    \end{align*}

We have now 7 types of terms that we again number from left to write and top to bottom. Again we collect terms whose first factor is of dimension $k.$ If $K=-1$ we have only the term number 2,
\begin{align*}
    \emptyset \otimes s\phi'_{2m+2} = s(\emptyset \otimes \phi'_{2m+2}).
\end{align*}

If $k=0$ we have the term number 1 and the term number 5 when $l=0$. Using Lemma~\ref{l:claim3} these two add up to
\begin{align*}
s\phi_0 \otimes \phi'_{2m+1} - \sum_{q=0}^{p-2}(s-s_q)\phi_{0} \otimes t^{q}d \phi'_{2m+1} = s(\phi_0 \otimes \phi_{2m+1}).  
\end{align*}

If $k>0$ is even then we have the term number 5 with $l=k/2$ and the term number 6 with $l=(k-2)/2$ which add up to (using Lemma~\ref{l:claim3} again)
\begin{align*}
    &s\phi_k\otimes t^{-(k-2)/2-1} \phi'_{2m-2(k/2-1)-1}- \sum_{q=0}^{p-2}(s-s_q)\phi_{k} \otimes t^{q-k/2}d \phi'_{2m-k+1}\\  
    &= s\phi_k\otimes t^{-k/2} \phi'_{2m-k+1} - \sum_{q=0}^{p-2}(s-s_q)\phi_{k} \otimes t^{q-k/2}d \phi'_{2m-k+1}\\
    &= s\left(\phi_k \otimes t^{-k/2}\phi'_{2m-k+1}\right).
\end{align*}

If $K$ is odd then we have the term number 3, when $l=(k-1)/2$, which is in the required form, and the terms number 4, when $l=(k-1)/2$ and number 7 when again $l=(k-1)/2$. The two latter terms cancel each other. Now by comparing the terms we have obtained for the three cases of $k$ with the formulas claimed by Proposition \ref{p:joinres} we see that if the formula for the even dimension claimed by the proposition is correct then the odd dimensional formula is correct. This finishes the inductive step of our proof of the proposition.

\section{Discussion}\label{s:discussion}

\subsection{Computer Implementation}
We have verified our main theorems on the Smith indices using computer calculations. We have written a code for computing the Smith classes in the mathematical programming language GAP~\cite{GAP4}, and utilizing the \texttt{simpcomp} package \cite{simpcomp}. 

We computed the Smith index for essentially the smallest complex where the top Smith class vanishes modulo 2. As described before, this is simply a $S^1$ with two capped Mobius bands (as in Figure~\ref{fig:mel}) glued  to it along the center of the Mobius bands. Let $K$ denote this $\Zring_2$-complex. We have verified that $K$ has non-zero integer Smith class of dimension 2, and this class is zero when computed modulo 2. Therefore, $I(K)=3$ and $I_2(K)=2$. We then computed the Smith index of $K*K$ and verified that the 5-dimensional Smith class is zero and the 4-dimensional class is non-zero, however is zero modulo~2. Thus $I(K*K)=I(K)+I(K)-1=5$. The three dimensional class is non-zero modulo 2, and thus $I_2(K*K) = I_2(K)+I_2(K)=4$.

We verified instability also for the complex where two capped 2-towers of Mobius bands are attached to the circle. For the join of this complex too the Smith class of dimension 5 vanishes, the class of dimension 4 is non-zero, but zero modulo 2, and the class of dimension 3 is non-zero modulo 2. It is interesting to note that, for the join of this simple complex with itself, the cochain computed by the program whose coboundary is equal to the Smith class of dimension 5, has integer coordinates with more than 100000 digits in a single coordinate and requires roughly 10GB of memory for the computation. As we raise the height of the tower of Mobius bands, the integer system becomes harder to solve.

We also verified that the 5-dimensional Smith class of the deleted join of the 2-dimensional Melikhov complex (depicted partly in Figure~\ref{fig:mel}) indeed vanishes modulo 2, but not as an integer class. However, computing the top Smith class of the join of two such deleted joins seems computationally infeasible due to huge space requirements.

\begin{figure}
    \centering
    \includegraphics[scale=0.7]{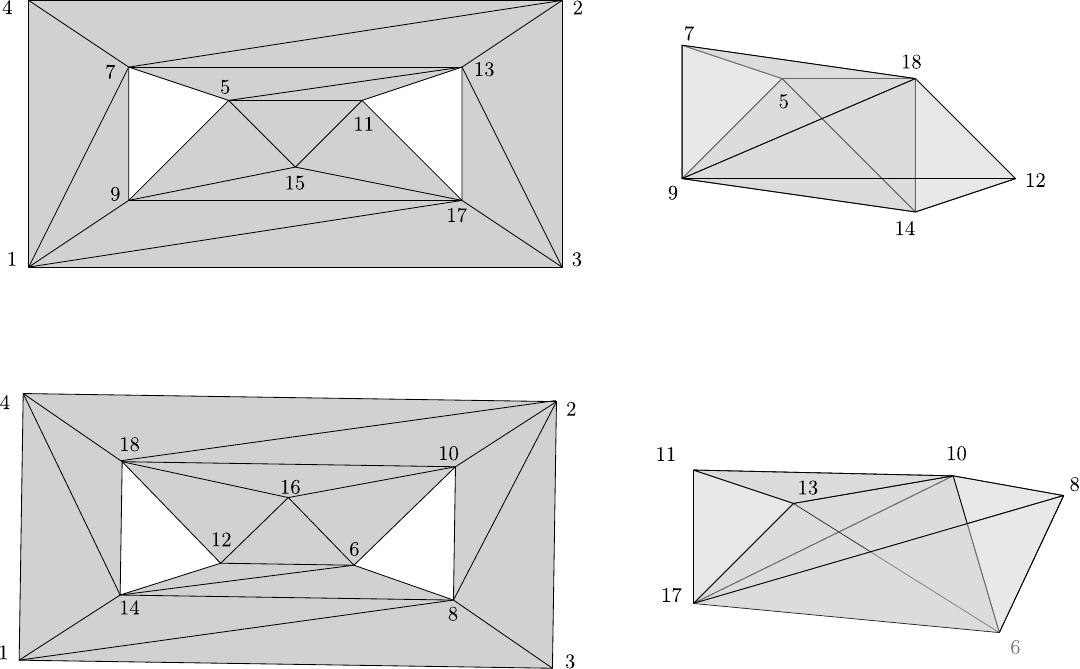}
    \caption{Complex of Example \ref{ex:b} is obtained by identifying vertices with the same label. The action sends an odd $j$ to $j+1$ and vice versa.}
    \label{fig:exb}
\end{figure}

The explicit simplicial complex of Example \ref{ex:b} which we have used in our computations is depicted in Figure~\ref{fig:exb}. The vertices with the same label are identified. As mentioned before, our computations show that the 2-dimensional Smith class of this $\Zring_2$-complex is non-zero modulo 2, however, the join of it with the complex of Example~\ref{ex:a} has vanishing 5-dimensional Smith class.

\subsection{Open Questions}
The most immediate, but a bit technical, open question is: Is it true that if the 5-dimensional Smith class of the deleted join of any 2-dimensional complex is non-zero modulo 2, then it has a boundary-equivariant dual mod $2^m$ for all $m$? In other words, the question asks if the instability can ever happen for the join of two deleted joins of 2-complexes, if at least one of the joins has non-zero top dimensional Smith class modulo $2$. 

Many simplicial complexes contain as a subcomplex a minimal complex in the sense of Sarkaria~\cite{Sar91b}. If this is the case, then the deleted product contains an antipodal sphere. It is reasonable that a complex which has non-zero modulo 2 van Kampen obstruction, and, does not contain a subcomplex which is minimal would provide an example where the above question is answered negatively. We conjecture that indeed this is the case and the answer to the above question is negative.

The most prominent open question related to this work is the generalized Menger embedding conjecture. Menger \cite{Men29} conjectured in 1929 that given two non-planar graphs $G, G'$ the product complex $G \times G'$ does not embed into the Euclidean 4-space. This conjecture was proved by Ummel \cite{Umm78} and later, using simpler arguments, by M. Skopenkov \cite{MSk03}. The natural generalization of this conjecture can be stated as:
\begin{Conjecture}
Let $K, L$ be complexes of positive dimensions, both with non-zero van Kampen obstruction. Then the van Kampen obstruction of $K \times L$ is non-zero.
\end{Conjecture}

See \cite{ASk14} for simple proofs and a survey of embeddability of products in low dimensions. Following the results of this paper it is possible that the conjecture is false, especially if there are geometric arguments relating the Smith index of $K \times L$ to that of $K*L$. However, we do not lean toward any side on this conjecture.

We have presented this paper using almost elementary algebraic topology, whereas more advanced language might shorten the arguments to some extent, for instance using spectral sequences. Moreover, we have not presented our results in their most generality, rather, we have restricted ourselves to simplicial complexes. These make the paper more accessible, however, leave generalizations of the results and the method for a future work. 

\section{Acknowledgments}
The work presented here has been done in a span of more than two years. In the first year the author was a post-doctoral researcher at the Institute for Research in Fundamental Sciences (IPM) in Tehran, Iran. In the second year, the author was a postdoctoral researcher of Computer Science at Saint Louis University, in Saint Louis, Missouri. I am grateful to these two institutions. I also thank my PhD advisor, Herbert Edelsbrunner, for igniting my interest in objects in 4-dimensional space and beyond. I also thank R. Karasev, F. Lazarus, and S. Melikhov for very helpful comments about the paper. I thank an anonymous referee for valuable suggestions which improved the presentation and the contents of the paper.

\bibliographystyle{plain}
\bibliography{bib.bib}

\begin{thebibliography}{10}

\bibitem{BKK02}
Mladen Bestvina, Michael Kapovich, and Bruce Kleiner.
\newblock Van {K}ampen’s embedding obstruction for discrete groups.
\newblock {\em Inventiones mathematicae}, 150(2):219--235, 2002.

\bibitem{Mat08}
A.~Bj{\"o}rner, G.M. Ziegler, and J.~Matou{\v{s}}ek.
\newblock {\em Using the Borsuk-Ulam Theorem: Lectures on Topological Methods
  in Combinatorics and Geometry}.
\newblock Universitext. Springer Berlin Heidelberg, 2008.

\bibitem{CoFl60}
P.~E. Conner and E.~E. Floyd.
\newblock Fixed point free involutions and equivariant maps.
\newblock {\em Bulletin of the American Mathematical Society}, 66(6):416--441,
  1960.

\bibitem{simpcomp}
F.~Effenberger and J.~Spreer.
\newblock {simpcomp}, a gap toolbox for simplicial complexes, {V}ersion 2.1.10.
\newblock \url {https://simpcomp-team.github.io/simpcomp/}
  {\texttt{https://simpcomp-team.github.io/}\discretionary
  {}{}{}\texttt{simpcomp/}}, Jun 2019.
\newblock Refereed GAP package.

\bibitem{Flo32}
R~Flores.
\newblock {\"U}ber die {E}xistenz {$n$}-dimensionaler {K}omplexe, die nicht in
  den {$R_{2n}$} topologisch einbettbar sind.
\newblock In {\em Ergebnisse eines Mathematischen Kolloquiums}, volume~5, pages
  17--24, 1932-33.

\bibitem{Freetal94}
Michael~H Freedman, Vyacheslav~S Krushkal, and Peter Teichner.
\newblock Van {K}ampen’s embedding obstruction is incomplete for 2-complexes
  in {$\mathbb{R}^4$}.
\newblock {\em Math. Res. Lett}, 1(2):167--176, 1994.

\bibitem{GAP4}
The GAP~Group.
\newblock {\em {GAP -- Groups, Algorithms, and Programming, Version 4.11.0}},
  2020.

\bibitem{GoGr2021}
Rafael Gomes and Gustavo Granja.
\newblock The cohomological index of free $\mathbb{Z}/p$-actions is not
  additive with respect to join.
\newblock {\em Topology and its Applications}, 291:107612, 2021.

\bibitem{Gru69}
Branko Gr{\"u}nbaum.
\newblock Imbeddings of simplicial complexes.
\newblock {\em Commentarii Mathematici Helvetici}, 44(1):502--513, 1969.

\bibitem{Hat02}
A.~Hatcher.
\newblock {\em Algebraic Topology}.
\newblock Algebraic Topology. Cambridge University Press, 2002.

\bibitem{Matetal10}
Ji{\v{r}}{\'\i} Matou{\v{s}}ek, Martin Tancer, and Uli Wagner.
\newblock Hardness of embedding simplicial complexes in {$\mathbb{R}^d$}.
\newblock {\em Journal of the European Mathematical Society}, 13(2):259--295,
  2010.

\bibitem{Mel09}
Sergey~A. Melikhov.
\newblock The van {K}ampen obstruction and its relatives.
\newblock {\em Proceedings of the Steklov Institute of Mathematics},
  266(1):142--176, Sep 2009.

\bibitem{MeSh06}
Sergey~A. {Melikhov} and Evgenij~V. {Shchepin}.
\newblock {The telescope approach to embeddability of compacta}.
\newblock {\em arXiv Mathematics e-prints}, December 2006.
\newblock math/0612085.

\bibitem{Men29}
Karl Menger.
\newblock {\"U}ber {P}l\"attbare {D}reiergraphen und {P}otenzen nicht
  pl\"attbarer {G}raphen.
\newblock In {\em Ergebnisse eines Mathematischen Kolloquiums}, volume~2, pages
  30--31, 1929.

\bibitem{Mil56}
John Milnor.
\newblock Construction of universal bundles, {II}.
\newblock {\em Annals of Mathematics}, 63(3):430--436, 1956.

\bibitem{Nak56}
Minoru Nakaoka.
\newblock Cohomology theory of a complex with a transformation of prime period
  and its applications.
\newblock {\em Journal of the Institute of Polytechnics, Osaka City University.
  Series A: Mathematics}, 7(1-2):51--102, 1956.

\bibitem{Pa20a}
Salman {Parsa}.
\newblock {On the Smith classes, the van Kampen obstruction and embeddability
  of $[3]*K$}.
\newblock {\em arXiv e-prints}, January 2020.
\newblock arXiv:2001.06478.

\bibitem{PaSk20}
Salman {Parsa} and Arkadiy {Skopenkov}.
\newblock {On embeddability of joins and their `factors'}.
\newblock {\em arXiv e-prints}, March 2020.
\newblock arXiv:2003.12285.

\bibitem{Sar91b}
Karanbir~S Sarkaria.
\newblock Kuratowski complexes.
\newblock {\em Topology}, 30(1):67--76, 1991.

\bibitem{Sch93}
G{\"o}ran Schild.
\newblock Some minimal nonembeddable complexes.
\newblock {\em Topology and its Applications}, 53(2):177--185, 1993.

\bibitem{Sha57}
Arnold Shapiro.
\newblock Obstructions to the {I}mbedding of a {C}omplex in a {E}uclidean
  space: {I}. {T}he first {O}bstruction.
\newblock {\em Annals of Mathematics}, pages 256--269, 1957.

\bibitem{Sko02}
Arkadiy Skopenkov.
\newblock On the {H}aefliger-{H}irsch-{W}u invariants for embeddings and
  immersions.
\newblock {\em Commentarii Mathematici Helvetici}, 77(1):78--124, 2002.

\bibitem{ASk14}
Arkadiy Skopenkov.
\newblock Realizability of hypergraphs and {R}amsey link theory.
\newblock {\em arXiv preprint arXiv:1402.0658}, 2014.

\bibitem{MSk03}
Mikhail Skopenkov.
\newblock Embedding products of graphs into {E}uclidean spaces.
\newblock {\em Fundamenta Mathematicae}, 179(3):191--198, 2003.

\bibitem{Spa89}
Edwin~H Spanier.
\newblock {\em Algebraic topology}.
\newblock Springer Science \& Business Media, 1989.

\bibitem{Ste51}
Norman Steenrod.
\newblock {\em The Topology of Fibre Bundles. (PMS-14)}.
\newblock Princeton University Press, 1951.

\bibitem{Umm73}
Brian~R Ummel.
\newblock Imbedding classes and {$n$}-minimal complexes.
\newblock {\em Proceedings of the American Mathematical Society},
  38(1):201--206, 1973.

\bibitem{Umm78}
Brian~R Ummel.
\newblock The product of nonplanar complexes does not imbed in 4-space.
\newblock {\em Transactions of the American Mathematical Society}, pages
  319--328, 1978.

\bibitem{vKam33}
Egbert~R Van~Kampen.
\newblock Komplexe in euklidischen {R}{\"a}umen.
\newblock In {\em Abhandlungen aus dem Mathematischen Seminar der
  Universit{\"a}t Hamburg}, volume 9-1, pages 72--78. Springer, 1933.

\bibitem{Web67}
Claude Weber.
\newblock Plongements de polyedres dans le domaine metastable.
\newblock {\em Commentarii Mathematici Helvetici}, 42(1):1--27, 1967.

\bibitem{Wu74}
W.~Wu.
\newblock {\em A {T}heory of {I}mbedding, {I}mmersion, and {I}sotopy of
  {P}olytopes in a Euclidean Space}.
\newblock Science Press, 1974.

\end{thebibliography}

\end{document}